\documentclass[a4paper,12pt]{article}
\usepackage[cp1251]{inputenc}
\usepackage{amssymb,amsmath,amsthm}
\usepackage[russian]{babel}
\usepackage{amssymb,amsmath,amsthm}
\usepackage{indentfirst}
\usepackage{hyperref,hypbmsec}
\usepackage[pdftex]{graphicx}
\numberwithin{equation}{section}
\newtheorem{Th}{\hskip\parindent Теорема}[section]
\newtheorem{Le}{\hskip\parindent Лемма}[section]

\newtheorem{Zam}{\hskip\parindent Замечание}[section]
\newtheorem{Hyp}{\hskip\parindent Гипотеза}[section]
\newtheorem{Op}{\hskip\parindent Определение}[section]

\newcommand{\A}{\mathcal{A}}

\newcommand{\R}{\mathfrak{R}}

\newcommand{\D}{\mathfrak{D}}
\newcommand{\N}{\mathbb{N}}
\newcommand{\M}{\mathfrak{M}}
\newcommand{\NN}{\mathfrak{N}}

\newcommand{\q}{\mathbf{q}}

\newcommand{\p}{\mathbf{p}}

\newcounter{propet}

\renewcommand{\le}{\leqslant}\renewcommand{\ge}{\geqslant}

\begin{document}

\author{И.\, Д.\, Кан 
(I.\,D.\,Kan)\footnote{Работа выполнена при поддержке РФФИ (грант  15-01-05700 А)} }

\title{
\begin{flushright}
\small{ Посвящается светлой памяти\\ Н. М. Коробова.}
\end{flushright}
\begin{flushleft}
УДК 511.321 + 511.31
\end{flushleft}Усиление \, \, теоремы Бургейна --- Конторовича \,-\, IV
(A strengthening of a theorem of Bourgain-Kontorovich- IV)}
\date{}
\maketitle
\begin{abstract}
Zaremba's conjecture (1971) states that every positive integer number $d$ can be represented as a denominator  of a finite continued fraction  $\frac{b}{d}=[d_1,d_2,\ldots,d_{k}],$ with all partial quotients $d_1,d_2,\ldots,d_{k}$ being bounded by an absolute constant $A.$ Several new theorems concerning this conjecture were proved by Bourgain and Kontorovich in 2011. The easiest of them states that the set of numbers satisfying Zaremba's conjecture with $A=50$ has positive proportion in $\N.$ In 2014 I. D. Kan and   D. A. Frolenkov proved this result with $A=5.$ In this paper the same theorem is proved with $A=4.$
 \noindent

\noindent
\textbf{Bibliography:} 13 titles. \par

\end{abstract}

\date{}
\maketitle
\begin{abstract}
В настоящей работе доказывается, что почти все натуральные числа являются знаменателями тех конечных цепных дробей,  все неполные частные которых принадлежат алфавиту $\{1,2,3,4\}$. Ранее аналогичная теорема была известна лишь для алфавитов большей мощности. Именно, впервые результат такого рода для алфавита 
$\{1,2,\ldots,50\}$ получили в 2011 году Бургейн и Конторович. Далее, в 2013 году автор статьи совместно с Д. А. Фроленковым доказали теорему для алфавита $\{1,2,3,4,5\}$.   Результат автора 2014 года, предшествующий настоящему, относился к алфавиту 
$\{1,2,3,4,10\}.$

 \noindent
\textbf
{Библиография:} 13 названий.\\
\noindent
\textbf
{Ключевые слова и выражения:} цепная дробь,  тригонометрическая сумма, гипотеза Зарембы. \par

\end{abstract}

\setcounter{Zam}0
\section{История вопроса}
Через $\left     [d_1,d_2,\ldots,d_{k}\right]$ обозначена конечная цепная дробь

\begin{equation}
\label{cont.fraction}
[d_1,d_2,\ldots,d_k]=\cfrac{1}{d_1+{\atop\ddots\,\displaystyle{+\cfrac{1}{d_k}}}}
\end{equation}
с натуральными неполными частными 
$\ d_1,d_2,\ldots,d_{k}$ (где $k$ --- натуральное), а через $\R_{\A}$ --- множество рациональных чисел $\frac{b}{d}$ , представимых  конечными цепными дробями  с неполными частными из некоторого конечного алфавита $\A \subseteq \N$:

  $$
\R_{\A}=\left\{\frac{b}{d}=[d_1,d_2,\ldots,d_{k}]\Bigl|{} \  d_j\in\A{}\ \mbox{для}{}\  j=1,\ldots,k\right\}.
$$
Через $\D_{\A}$    обозначено множество знаменателей $d$  чисел $\frac{b}{d}\in\R_{\A}$, а через $\D_{\A}(N)$ --- множество таких знаменателей, ограниченных сверху числом $N\in \N$: 

$$\D_{\A}=\left\{d\in \N\Bigl|{} \  \exists b:{}\  \gcd(b,d)=1, {}\  \frac{b}{d}\in\R_{\A}\right\},{}\ {}\ {}\ {}\ \D_{\A}(N)=\left\{d\in \D_{\A}\Bigl|{} \   d\le N\right\}.$$

\begin{Hyp}\label{h1.1} (Гипотеза Зарембы  \cite{Zaremba}). Существует константа $A$ (скорее всего, $A=~5$), такая что для любого $N\in \N$ для алфавита 
\begin{equation}
\label{AA}
\A= 1,2,...,A 
\end{equation}
имеет место равенство 
$|\D_{\A}(N)|=N$.\end{Hyp}

Обзор результатов, связанных с гипотезой 1.1, можно найти в  работах \cite{BK},\cite{NG}. 

Для фиксированного алфавита $\A$ число $d$ называется допустимым \cite{BK}, если  для любого $q>1$  множество $\D_{\A}$ содержит число, сравнимое с $d$ по модулю $q$. Множество допустимых чисел обозначаено через $\mathfrak{A}_{\A}$.
Пусть  $\delta_{\A}$ --- хаусдорфова размерность множества бесконечных цепных дробей с неполными частными из алфавита ${\A}$. Бургейн и Конторович в 2011 году доказали, в частности, следующее.

\begin{Th}\label{t1.1} \cite[ теорема 1.25]{BK}. Для каждого
алфавита ${\A}$, такого что
\begin{equation}
 \label{c}                                                                                                                                                                                                             
 \delta_{\A}>\frac{307}{312}=0.9839\ldots,                                                                                                                                   
\end{equation}
справедливо неравенство (``положительная пропорция''):
\begin{equation}
 \label{lc}                                                                                                                                                                                                             
 |\D_{\A}(N)|>> N.                                                                                                                                   
\end{equation}
\end{Th}

\begin{Th}\label{t1.2}
 \cite[ теорема 1.27]{BK}.  Для каждого алвавита ${\A}$, удовлетворяющего условию (\ref{c}), существует подмножество $\widetilde{\D}_{\A} \subseteq \D_{\A}$, содержащее почти все допустимые числа. То есть,  найдется константа $c=c({\A})>0$, такая что
 
 \begin{equation}
 \label{lld}                                                                                                                                                                                                             
\frac{\left|\widetilde{\D}_{\A}\bigcap\left [\frac{N}{2},N\right] \right|}
{\left|\mathfrak{A}_{\A}\bigcap\left [\frac{N}{2},N\right] \right|}
= 1+
 O\left(
 \exp{\left\{-c\sqrt{\log N}\right\}}
 \right).                                                                                                                              
 \end{equation}
Следовательно, равенство (\ref{lld}) имеет место при замене $\widetilde{\D}_{\A}$ на ${\D}_{\A}$. Кроме того, каждое число $d\in \widetilde{\D}_{\A}$  появляется  с кратностью
 \begin{equation}
 \label{cdd}                                                                                                                                                                                                             
>>{}\ N^ {2\delta_{\A}-{1.001}}.                                                                                                                                   
\end{equation}
\end{Th}
     
    Далее теорема \ref{t1.1} усиливалась в работах  \cite{FK1} --- \cite{FK4}, а теорема \ref{t1.2} --- в работах \cite{Huang} и \cite{K5}.
    (Конечно, теорема \ref{t1.1} следует из теоремы \ref{t1.2}, но допускает также и более простое доказательство.) В частности, в работе \cite{K5} было доказано, что в теоремах \ref{t1.1}  и \ref{t1.2} неравенство  (\ref{c}) можно заменить  более слабым условием $\delta_{\A}>0.8$, которому удовлетворяют все алфавиты вида  $\{1,2,3,4,n\}$ при $n$, принимающем любое из значений $6,\ldots,                                                                                                                                    10$.
    
    Некоторый условный результат по проблеме имеется также в работе \cite{Mag}.

\section{Основные результаты работы}

В  настоящей статье имеется две основных теоремы. 

\begin{Th}\label{2.1} Для произвольного алфавита $ \A$, такого что 
 \begin{equation}
 \label{hy}                                                                                                                                                                                                             \delta_{\A}>\frac{11}{14}=0.7857\ldots ,                                                                                                                                   
 \end{equation}
 справедливо неравенство $ |\D_{\A}(N)|>> N$.                                                                                                                                   
\end{Th}

\begin{Th}\label{2.2}
 Пусть алфавит $\A$ удовлетворяет неравенству (\ref{hy}). Тогда 
 существует константа $c=c({\A})>0$, такая что имеют место формулы (\ref{lld})                                                                                                                                                                                                           и (\ref{cdd}).
\end{Th}

\begin{Zam}\label{z2.1}
  Согласно результатам  Дженкинсона \cite{Jenkinson}, неравенству (\ref{hy}) удовлетворяет  алфавит (\ref{AA}) при $A=4$.
  \end{Zam}

\section{Обозначения} 
Всюду далее $\varepsilon_0~\in~\left  (0,{}\ 0.0004\right)$  --- произвольно малая положительная константа, участвующая в построении ансамбля $\Omega^{(N)}=\Omega^{(N,\varepsilon_0)}$. Знак Виноградова $f(N)\ll g(N)$ для двух произвольных функций $f(N)$ и $g(N)$ обозначает существование константы  $C,$ зависящей только от $\A$ и $\varepsilon_0,$ такой что  $|f(N)|\le Cg(N).$   Также используются традиционные обозначения  $e(x)=\exp(2\pi ix),{}\ e_n(x)=\exp\left(\frac{2}{n}\pi ix\right).$ Наибольший общий делитель двух целых чисел $a$ и $b$ обозначается через $\gcd(a,b).$ Для натурального $n$ и целого $m$
\begin{equation}
 \label{11hy}                                                                                                                                                                                                                  \delta_n(m)=
\frac{1}{n}\sum\limits^n_{k=1}e_n(km)=
\left\{
              \begin{array}{ll}
1,{}\ {}\ \mbox{если}{}\ {}\ m {}\ {}\ \mbox{делится на}{}\ {}\ n,                \\
               0,{}\ {}\ \mbox{--- в противном случае,}                                                                      
              \end{array}
\right.
       \end{equation}
--- $\delta$-символ Коробова (в честь Н. М. Коробова, пропагандировавшего идею использования формулы (\ref{11hy}), см.  \cite{Korobov}).  
   Мощность конечного множества  $S$ обозначается через $|S|$. Для действительного числа $\alpha$ через  $[\alpha]$,$\{\alpha\}$ и $\|\alpha\|$ обозначаются, соответственно,  целая часть от  $\alpha$, дробная доля $\alpha$ и расстояние от $\alpha$ до ближайшего целого:
$$[\alpha]=\max\left\{z\in\mathbb{Z}|\,z\le \alpha \right\},{}\ {}\   \{\alpha\}=\alpha-[\alpha],{}\ {}\ 
\|\alpha\|=\min\left\{\{\alpha\},\{-\alpha\}\right\}. $$
Кроме того, если $g$ \--- матрица, то $||g||$ \--- ее норма (определенная ниже в параграфе \ref{6}).

\section{Благодарности} Автор благодарит профессора Н. Г. Мощевитина за постановку задачи и неоднократное обсуждение темы статьи. Также автор благодарен Д. А. Фроленкову за многократное обсуждение и многие  полезные советы.

\section{Основные свойства ансамбля $\Omega^{(N)}$}\label{6}
Через $\Gamma_{\A}$ обозначена мультипликативная полугруппа $\Gamma_{\A}\subseteq SL\left(2,\mathbb{Z}\right)$ с единицей  $E=
\begin{pmatrix}
1 & 0 \\
0 & 1
\end{pmatrix}$, порожденная  матрицами 
$
\begin{pmatrix}
1   & v \\
u & uv+1
\end{pmatrix},
$
где $u,  v \in\A.$  
Для  четного $k$ и чисел $d_1,d_2,\ldots,d_{k}\in\A$ в качестве нормы произвольной матрицы $g=\left( a{}\ {}\ b \atop{c{}\ {}\ d}\right)\in \Gamma_{\A}$, такой что 
$$
\begin{pmatrix}
a & b \\
c & d
\end{pmatrix}=
\begin{pmatrix}
1 & d_2 \\
d_1 & d_1d_2+1
\end{pmatrix}
\begin{pmatrix}
1 & d_4 \\
d_3 & d_3d_4+1
\end{pmatrix}\ldots
\begin{pmatrix}
1 & d_k \\
d_{k-1} & d_{k-1}d_{k}+1
\end{pmatrix},
$$
 рассматривается, как обычно (\cite{FK3} --- \cite{K5}),  величина
$||g||=d=\langle d_1,d_2,\ldots,d_{k}\rangle
$
 ---  знаменатель цепной дроби (\ref{cont.fraction}), не сократимый с ее числителем.

Скажем, что для некоторого множества $\Omega \subseteq \Gamma_{\A}$ имеет место разложение 
 
\begin{equation}
\label{1morm}
 \Omega=\Omega_1\Omega_2\Omega_3\ldots\Omega_n
 \end{equation}
на независимые множители   $ \Omega_1,\Omega_2,\Omega_3,\ldots,\Omega_n\subseteq \Gamma_{\A}$, если для каждой матрицы $g\in\Omega$ найдется, причем единственный, набор матриц $g_1,g_2,g_3,\ldots,g_n,$ таких что  $g_i\in \Omega_i$ для $i=1,2,3,\ldots,n$ и выполнено равенство
$$g=g_1g_2g_3\ldots g_n.
$$
Конечно, при этом выполняется равенство  
$ \left|\Omega\right|=|\Omega_1||\Omega_2||\Omega_3|\ldots|\Omega_n|.
 $

Всюду далее будем использовать обозначения $A=\max\A$,
$Q_1= \left[\exp{ \left  (A^4\varepsilon_0^{-5}\right)}\right]+1 .
 $    Положим $Q_0=0$    и определим последовательность $\{Q_j\}$ для $j$ от нуля до бесконечности:
 \begin{equation} 
\label{24norm}
  \left  \{Q_j\right\}^\infty_{j=0}= \left  \{0,Q_1,Q^2_1,Q^3_1,\ldots,Q^j_1,\ldots\right\}.
 \end{equation}
 Рассмотрим  две  произвольных  матрицы
$  g_2\in\Omega_2,{}\  {}\  g_4\in\Omega_4$, три  параметра $M_1,M^{(2)},M^{(4)}\in \mathbb{R}_{+}$
 и  следующие два неравенства:

\begin{equation}
\label{6unorm}
 \frac{ M^{(2)} }{150A^2{\left  (M_1\right)}^{2\varepsilon_0}} \le||g_2||\le 73A^2   M^{(2)}{ \left  (M_1M^{(2)}\right)}^{2\varepsilon_0},
 \end{equation}

\begin{equation}
\label{7anorm}
\frac{{ \left  (M^{(4)}\right)}^{1-\varepsilon_0} }{150A^2} \le||g_4||\le 73A^2M^{(4)}.
 \end{equation}
По достаточно большому числу $N$ и по малому параметру $\varepsilon_0~\in~\left  (0,{}\ 0.0004\right)$  в   \cite{FK3} было построено специальное множество матриц --- ансамбль  (см. терминологию в \cite{BK}) $$\Omega^{(N)}=\Omega^{(N,\varepsilon_0)}\subseteq \left\{g\in\Gamma_{\A}\,\Bigl|\,||g||\le1,02 N\right\},$$
 для которого имеет место разложение на независимые множители  (\ref{1morm}) с $n=4$ со свойствами, перечисленными в следующей лемме.

  \begin{Le}\label{l6.1}
 \cite[теорема 6.1]{K5} Существует  непустое множество матриц --- ансамбль $\Omega^{(N)}\subseteq \Gamma_{\A}$, такое что для  всякого $M_1 \in  [Q_1,N]$ найдeтся разложениe 
  \begin{equation}
\label{1norm}
 \Omega^{(N)}=\Omega_1\Omega
 \end{equation}
   на независимые множители  $\Omega_1$ и $\Omega$, для которых выполнен ряд свойств: 
   
    во-первых, имеет место оценка 
\begin{equation}
\label{2norm} 
 \left   |\Omega_1\right|>>{ \left  (M_1\right)}^{2\delta_{\A}-\varepsilon_0},
  \end{equation}
 
 во-вторых, для любых двух матриц $g_1\in\Omega_1$ и 
$g\in\Omega$  выполняются неравенства 
 \begin{equation}
\label{normyy}
\frac{M_1}{70A^2} \le||g_1||\le 1.01(M_1)^{1+2\varepsilon_0},{}\ {}\ {}\ {}\ 
\frac{N}{160A^2(M_1)^{1+2\varepsilon_0}} \le||g||\le 73A^2\frac{N}{M_1},
 \end{equation}

    в-третьих, для произвольных действительных чисел $
  M^{(2)},M^{(4)}\in  [Q_1,N]\bigcup\{1\}  $,  удовлетворяющих неравенству
\begin{equation}
\label{yorm}
M_1M^{(2)}M^{(4)}\le N,
  \end{equation}
  найдется разложение множества $\Omega$ вида
  \begin{equation}
\label{1niorm}
 \Omega=\Omega_2\Omega_3\Omega_4
 \end{equation}
   на независимые множители   $\Omega_2, \Omega_3, \Omega_4$,  для которых  выполнены как неравенства
   \begin{equation}
\label{2n1orm}
 \left   |\Omega_2\right|>>
 { \left  (M^{(2)}\right)}^{2\delta_{\A}} 
 { \left  (M_1M^{(2)}\right)}^{-2\varepsilon_0},{}\ 
 \left  |\Omega_4\right |>>{ \left  (M^{(4)}\right)}^{2\delta_{\A}-2\varepsilon_0},
  \end{equation}  
так и оценки (\ref{6unorm}), (\ref{7anorm})  --- для любых двух  матриц $  g_2\in\Omega_2,{}\  {}\  g_4\in\Omega_4$.
 В частности, если $M^{(2)}=1$ или $M^{(4)}~=~1,$ то $\Omega_2=\{E\}$ или $\Omega_4=\{E\}$, соответственно.
\end{Le}

Пусть число  $M_1 \in  [Q_1,N]$ уже как-либо выбрано, так что имеет место разложение (\ref{1norm}) со свойствами (\ref{2norm}) и (\ref{normyy}). В этом случае множество $\Omega$ из (\ref{1norm}) будем называть полуансамблем.
Если имеет место разложение  (\ref{1niorm}), то для любых двух элементов $g^{(1)}$ и $g^{(2)}$  полуансамбля $ \Omega$ введем обозначения 
$$g^{(1)}=g^{(1)}_2g^{(1)}_3g^{(1)}_4,{}\ {}\ {}\ 
g^{(2)}=g^{(2)}_2g^{(2)}_3g^{(2)}_4,
$$
  где нижний индекс $i {}\ (i=2,3,4)$ указывает на принадледность соответствующему множеству $\Omega_i$.  Далее,  если $X$ --- некоторое множество $2\times2$-матриц $g$, то $\widetilde{X}$ --- множество  вектор-столбцов $\widetilde{g}=g\left(0\atop{1}  \right)
$.        
 Для координат трех пар  произвольных векторов
  $$\widetilde{g}^{(1)},\widetilde{g}^{(2)}\in \widetilde{\Omega},{}\ {}\ {}\ {}\ 
  \widetilde{g}^{(1)}_2,\widetilde{g}^{(2)}_2\in \widetilde{\Omega}_2,{}\ {}\ {}\ {}\ 
  \widetilde{g}^{(1)}_4,\widetilde{g}^{(2)}_4\in \widetilde{\Omega}_4
  $$  
   введем такие обозначения:
\begin{equation}
\label{11norm}
 \widetilde{g}^{(1)}=\begin{pmatrix}
 x \\
X
\end{pmatrix},{}\ 
\widetilde{g}^{(2)}=\begin{pmatrix}
y \\
Y
\end{pmatrix},{}\ {}\ 
\widetilde{g}^{(1)}_2=\begin{pmatrix}
 x_2  \\
X_2
\end{pmatrix},{}\ 
\widetilde{g}^{(2)}_2=\begin{pmatrix}
y_2 \\
Y_2
\end{pmatrix},{}\ {}\ 
\widetilde{g}^{(1)}_4=\begin{pmatrix}
x_4 \\
X_4
\end{pmatrix},{}\ {}\ 
\widetilde{g}^{(2)}_4=\begin{pmatrix}
y_4 \\
Y_4
\end{pmatrix}.
\end{equation}
  
  Пусть $M_2,M_4\in \mathbb{R}_+$. Рассмотрим неравенства 
  \begin{equation}
\label{6unorm1}
 \frac{ M_2
 { \left  (M_1M_2\right)}^{-6\varepsilon_0} }{Q_5} \le ||g_2||\le 
 {M_2},
 {}\ {}\ {}\ {}\ {}\ {}\ 
 \frac{{ \left  (M_4\right)}^{1-\varepsilon_0} }{Q_5} \le||g_4||\le {M_4}.
 \end{equation}  
 Свойства полуансамбля $\Omega$ из леммы \ref{l6.1} несколько уточняет  следующая

 \begin{Th}\label{t6.1}
  Для любого числа $M_1
 \in  [Q_1,N]$ найдeтся   непустое множество матриц --- полуансамбль $\Omega\subseteq \Gamma_{\A}$, такое что   для произвольного числа $
  M_2\ge 1   $,  удовлетворяющего неравенству
$M_1M_2\le Q_3 N,
 $ найдется разложение  вида $\Omega=\Omega_2 \Omega_{3,4}$   
   на независимые множители   $\Omega_2$ и $\Omega_{3,4}$,  для которых выполнен ряд свойств: 
   
    (i) имеет место  неравенство
   \begin{equation}
\label{2n1orm1}
 \left   |\Omega_2\right|>>
 { \left  (M_2\right)}^{2\delta_{\A}} 
 { \left  (M_1M_2\right)}^{-10\varepsilon_0},
  \end{equation}  
  
(ii) для любой  матрицы $g_2\in\Omega_2$ имеет место первая из  оценок (\ref{6unorm1}):  в частности, в обозначениях (\ref{11norm}) выполнено неравенство 
$\max\left \{X_2,Y_2\right\}\le {M_2},
 $
 
    (iii) для любого числа $
  M_4  \ge 1  $,  удовлетворяющего неравенству
\begin{equation}
\label{1yorm1}
M_1M_2M_4\le Q_6 N,
  \end{equation}
    найдется разложение множества $\Omega_{3,4}$ вида $\Omega_{3,4}=\Omega_3 \Omega_4$   
   на независимые множители   $\Omega_3$ и $\Omega_4$,  такие что  выполнены  как  неравенство
  \begin{equation}
\label{2n1orm2}
\left  |\Omega_4\right |>>{ \left  (M_4\right)}^{2\delta_{\A}-2\varepsilon_0},
    \end{equation}
так  и вторая из оценок  (\ref{6unorm1})  --- для любой  матрицы $g_4\in\Omega_4$: в частности,  
 \begin{equation}
\label{normyy1}
\max\left \{X_4,Y_4\right\}\le {M_4}.
 \end{equation}
 \end{Th} 
 Доказательство. 
 Пусть выбраны значения величин $M_2$ и    $M_4$, удовлетворяющих условиям теоремы.  Тогда, используя обозначения  (\ref{24norm}),  положим:
$$\mathcal{M}^{(2)}=\frac{M_2}{Q_3}\left(M_1M_2\right)^{-4\varepsilon_0},{}\ {}\ {}\ {}\ 
\mathcal{M}^{(4)}=\frac{M_4}{Q_3}.
$$
Участвующие в лемме \ref{l6.1} величины $M^{(2)}$ и $M^{(4)}$ определим правилами:

$
M^{(2)}=
                {\mathcal{M}}^{(2)},$ если  ${\mathcal{M}}^{(2)}\ge Q_1,$                
   и           $M^{(2)}=  1$--- в противном случае;  аналогично,

             $M^{(4)}=
               \mathcal{M}^{(4)},$ если  $\mathcal{M}^{(4)}\ge Q_1,$
         и    $M^{(4)}=  1$ --- в противном случае.
          
    \noindent      Тогда числа  $M^{(2)}$ и $M^{(4)}$ принадлежат множеству $ [Q_1,N]\bigcup\{1\}$. Кроме того,  если ${M}^{(2)}=~1$ или ${M}^{(4)}=1$, то, полагая $\Omega_2=\{E\}$ или $\Omega_4=\{E\}$, соответственно, получаем, что все требуемые неравенства выполнены. Если же равенства ${M}^{(2)}=1$ или ${M}^{(4)}=1$ не выполнены, то, следовательно, для чисел  ${M}^{(2)}$ или ${M}^{(4)}$ справедливы равенства $
M^{(2)}=
                {\mathcal{M}}^{(2)}$ или $
M^{(4)}=
                {\mathcal{M}}^{(4)}$. Для таких значений 
   условие  (\ref{yorm}) леммы \ref{l6.1}  выполнено ввиду неравенства (\ref{1yorm1}). Поэтому доказаны  оценки (\ref{6unorm}), (\ref{7anorm}) и (\ref{2n1orm}).
Подставляя в них   значения ${M}^{(2)}$ и ${M}^{(4)}$, получаем неравенства (\ref{6unorm1}), (\ref{2n1orm1}) и (\ref{2n1orm2}). Теорема доказана.

Отметим, что далее довольно часто в качестве значений ${M}_2$  и ${M}_4$ выбираются числа  
\begin{equation}
\label{17narm}
  {M}_2=1,{}\ {}\ {}\ {}\    M_4=\max\left\{1, \frac{1}{2}Q_{\alpha-1}\right\} ,
 \end{equation}
где $\alpha\in \N$. В этом случае для проверки оценки (\ref{1yorm1}) достаточно установить неравенство
\begin{equation}
\label{17narb}
  {M}_1 Q_{\alpha} \le 
 {Q_6N}
 , 
 \end{equation}
 гарантирующее выполнение условий теоремы для числа ${M}_1\in  [Q_1,N]$.

\section{Основа доказательства формул (\ref{lc})~---~(\ref{cdd}).}\label{7}
Напомним обозначения из \cite{K5}. Применяя теорему Дирихле \cite[лемма 2.1, стр. 17]{Von}, для каждого $\Theta\in[0,1)$ найдем  целые числа $a$, $q$ и $l$  и действительное число $\lambda$, такие что 
\begin{equation}
\label{17norm}
  \Theta=\left\{\frac{a}{q}+\frac{l}{2N}+\frac{\lambda}{N}\right\}
  ,{}\ \gcd(a,q)=1, {}\ 0\le a<q\le\frac{\sqrt{N}}{Q_1},{}\ {}\ 
  \lambda\in \left (-\frac{1}{4},\frac{1}{4}\right ],{}\ {}\ 
  |l|\le\frac{3}{q}Q_1\sqrt{N},
 \end{equation}
 при чем равенство $a=0$ возможно только при $q=1$.   
 Фиксируем константу
 $T_1=7Q_7
$ (где $Q_7$ --- элемент последовательности (\ref{24norm})) и целое число $\kappa\in\left[0,{}\ T_1-1\right]$. Для целых $\alpha,\beta\ge 1$ рассмотрим числа $\Theta$ из (\ref{17norm}), удовлетворяющие соотношениям 
 \begin{equation}
\label{43norm}
 l\equiv\kappa\pmod {T_1},{}\ {}\ {}\ 
Q_ {\alpha-1}< q\le Q_{\alpha} , {}\ {}\ 
Q_ {\beta-1}\le         \left|l \right|         \le Q_{\beta} ,
  \end{equation}  
и положим   
$$P_{\alpha,\beta}=P^{(\lambda)}_{\alpha,\beta}(\kappa)=\left \{\Theta{}\ \Bigl|{} \ 
 \mbox{выполнены}{}\  (\ref{17norm}){}\  \mbox{и}{}\  (\ref{43norm})  \right\}.
$$

   Всюду далее   $Z$ --- произвольное непустое подмножество конечного множества $P_{\alpha,\beta}$, которое предполагается  непустым.
   Положим также
$$\sigma_{N,Z}=\sum_{\Theta\in Z}\left|
\sum_{g\in\Omega^{(N)}}e\left((0,1)\widetilde{g}\Theta\right)
\right|.
$$
С помощью метода Хуанга \cite{Huang} (обобщившего методы  Бургейна --- Конторовича  \cite{BK} и других \cite{FK4}) в \cite{K5} была доказана следующая

\begin{Le}\label{l7.1}(\cite[теорема 7.1]{K5}) Пусть найдется не зависящая от $\varepsilon_0$ константа $c=c(\A)>0$, такая что выполняется оценка  
\begin{equation} 
\label{45norm}
\sigma_{N,Z}\ll
   \frac{\left   |\Omega^{(N)}\right|\sqrt{|Z|}}
  {\left(Q_{\alpha}Q_{\beta}\right)^{c+O(\varepsilon_0)}}.
\end{equation}
 Тогда для алфавита $\A$ имеют место формулы (\ref{lc})~---~(\ref{cdd}).
 \end{Le}

Для двух произвольных  чисел $\Theta^{(1)},\Theta^{(2)}\in Z\subseteq P_{\alpha,\beta}$ введем обозначения 
\begin{equation} 
\label{50norm}
\Theta^{(1)}=\frac{a^{(1)}}{q^{(1)}}+\frac{l^{(1)}}{2N}+\frac{\lambda}{N}, {}\ {}\  {}\ {}\ \Theta^{(2)}=\frac{a^{(2)}}{q^{(2)}}+\frac{l^{(2)}}{2N}+\frac{\lambda}{N}
\end{equation}
и положим
\begin{equation} 
\label{67mnorm}
\p=\gcd(q^{(1)},q^{(2)}),{}\ {}\ {}\ {}\ \q=\frac{1}{\p}q^{(1)}q^{(2)}.
\end{equation}
Следовательно, полагая
$q^{(1)}_0=\frac{1}{\p}q^{(1)},{}\ {}\ {}\ q^{(2)}_0=\frac{1}{\p}q^{(2)},
$ получим
$\q=\p q^{(1)}_0q^{(2)}_0.
$

Далее, напомним обозначения (\ref{11norm}) и обозначим  через $t$  и $T$ числители дробей $$\left\|  x\frac{a^{(1)}}{q^{(1)}}-y\frac{a^{(2)}}{q^{(2)}}\right\|=
\frac{t}{\q},{}\ {}\ {}\ {}\ 
\left\|  X\frac{a^{(1)}}{q^{(1)}}-Y\frac{a^{(2)}}{q^{(2)}}\right\|=
\frac{T}{\q}.
$$
Другими словами,  
\begin{equation} 
\label{69n2f1orm}\begin{array}{ll}
\left|xa^{(1)}{q^{(2)}_0}-ya^{(2)}q^{(1)}_0\right|
\equiv t \pmod{\p q^{(1)}_0
q^{(2)}_0},{}\ {}\  0\le t<\q, \\{}\ \\
\left|Xa^{(1)}{q^{(2)}_0}-Ya^{(2)}q^{(1)}_0\right|\equiv T \pmod{\p q^{(1)}_0
q^{(2)}_0},{}\ {}\  0\le T<\q .
\end{array}\end{equation}
Рассмотрим неравенства
\begin{equation} 
\label{52norm}
\frac{t}{\q}\le\min
\left\{\frac{1}{M_1}74A^2Q_{\beta},{}\ {}\ {}\ 
                              \frac{1}{N}150A^3x+\left\|  \frac{1}{2N}\left(xl^{(1)}-yl^{(2)}\right)\right\|+
                \left\| \frac{\lambda}{N}\left(  x-y\right)\right\|    
            \right\},
\end{equation}

\begin{equation} 
\label{52norm1}
\frac{T}{\q}\le\min
\left\{ \frac{1}{M_1}74A^2Q_{\beta},{}\ {}\ {}\ 
                              \frac{1}{N}150A^3X+\left\|  \frac{1}{2N}\left(Xl^{(1)}-Yl^{(2)}\right)\right\|+
                \left\| \frac{\lambda}{N}\left(  X-Y\right)\right\|     
             \right\},
\end{equation}

\begin{equation} 
\label{53norm}
\left| xl^{(1)}-yl^{(2)} \right|\le\left(9A \right) ^5x+2N\frac{t}{\q},{}\ {}\ {}\ {}\ 
\left| Xl^{(1)}-Yl^{(2)} \right|\le\left(9A \right) ^5X+2N\frac{T}{\q}.
\end{equation}

Введем также обозначения 
$$\NN=\left\{    \left(\widetilde{g}^{(1)},\widetilde{g}^{(2)},\Theta^{(1)},\Theta^{(2)}\right) \in   \widetilde{\Omega}^2\times Z^2 \Bigl|{} \  \mbox{выполнены неравенства} {}\ \mbox{(\ref{50norm})}{}\  \mbox{---}{}\ (\ref{53norm})\right\},
$$ 
$$\NN(\p,t,T)=\left\{    \left(\widetilde{g}^{(1)},\widetilde{g}^{(2)},\Theta^{(1)},\Theta^{(2)}\right) \in   \NN \Bigl|{} \  \mbox{значения параметров}{}\ \p,t,T{}\  \mbox{фиксированы} \right\}.$$

\begin{Le}\label{l7.3}  (\cite[см. лемму 8.5]{K5}). Если выполнено неравенство
 \begin{equation} 
\label{69n1fo2r3m}
 75A^2 Q_{\alpha}Q_{\beta}{}\ {}\ \le {}\ {}\ M_1{}\ {}\ \le {}\ {}\ 
  \min\left\{\left(  Q_{\alpha}Q_{\beta}   \right)^5   ,{}\ 
  \frac{N}{\sqrt{Q_{\alpha}Q_{\beta}}}\right\},
\end{equation}
 то имеет место оценка
\begin{equation} 
\label{55norm}
\sigma_{N,Z}\ll \left(M_1\right)^{1+2\varepsilon_0}   
\sqrt{\left|\Omega_1\right|\left| \NN   \right| }
   .
\end{equation}\end{Le}

Доказательство. Для получения оценки (\ref{55norm}) из аналогичного неравенства работы \cite{K5} временно  положим 
\begin{equation} 
\label{54nowrm}
\Omega_2=\Omega_3=\{E\},{}\  
\Omega_4=\Omega,
 \end{equation}
 тогда  знак суммы по множеству $\Omega_3$ теперь не нужен, и неравенство (\ref{55norm}) при выполнении условий (\ref{54nowrm}) доказано.  
  Теперь покажем, что  от выбора параметров (\ref{54nowrm}) можно  отказаться. Действительно, для этого достаточно применить   теорему \ref{t6.1} еще раз, полагая $\Omega=\Omega_4$ и придавая обозначениям $\Omega_2, \Omega_3, \Omega_4$ другие значения. Лемма доказана.

Далее для краткости $\delta_{\A}$ обозначается через $\delta$. 

\begin{Th}\label{t7.1} Пусть  найдется такая не зависящая от $\varepsilon_0$ константа
$\mathbf{c}=\mathbf{c}(\A)>0$, при которой для любых натуральных значений $\alpha$ и $\beta$ существует число $M_1=M_1(\alpha,\beta)$ из интервала (\ref{69n1fo2r3m}),  такое что    
\begin{equation} 
\label{70nokrm}
\frac{\left| \NN \right|}{|Z||\Omega|^2}
\ll 
\left(M_1
 \right)^{-2+2\delta-\mathbf{c}+O(\varepsilon_0)}.
 \end{equation}
Тогда для алфавита $\A$ имеют место формулы (\ref{lc})~---~(\ref{cdd}). 
\end{Th}
Доказательство. Ввиду неравенств (\ref{2norm}) и (\ref{70nokrm}),  с помощью леммы \ref{l7.3}
легко получить оценку (\ref{45norm}). Теперь утверждение теоремы следует из леммы \ref{l7.1}. Теорема доказана.

Здесь и далее   $\mathbf{c}=\mathbf{c}(\A)>0$ --- произвольная достаточно малая   не зависящая  от $\varepsilon_0$ константа, о которой  будет доказано неравенство (\ref{70nokrm}). Можно, для определенности, считать, что  $\mathbf{c}=0.001\left(\delta-\frac{11}{14}\right),$ но это равенство нигде не будет использовано.

\section{Оценка величины $ \left| \NN \right|$ суммой мощностей множеств}

Напомним обозначения (\ref{11norm}), (\ref{50norm}) и  (\ref{69n2f1orm}) и обозначим:
$$\M=\left\{    \left(\widetilde{g}^{(1)},\widetilde{g}^{(2)},\Theta^{(1)},\Theta^{(2)}\right) \in  \NN  \Bigl|{} \  {}\   T=t=0 \right\}.
$$
Тогда,
 ввиду (\ref{53norm}), на всех элементах множества $\M$ выполнены неравенства

\begin{equation} 
\label{66norm}
\left| xl^{(1)}-yl^{(2)} \right|\le\left(9A \right) ^5x,{}\ {}\  
\left| Xl^{(1)}-Yl^{(2)} \right|\le\left(9A \right) ^5X.
\end{equation}

\begin{Zam}\label{z9.1}  В работе \cite[пп. 6.2.1]{BK} (см. также \cite[доказательство леммы 14.1]{FK3}) было показано, что условие $T=t=0$ равносильно следующему:
\begin{equation} 
\label{67norm}
q^{(1)}=q^{(2)}=\q=\p ,{}\ xa^{(1)}\equiv ya^{(2)}\pmod{\q},{}\ 
Xa^{(1)}\equiv Ya^{(2)}\pmod{\q}
.
\end{equation}
 Кроме того, при фиксированном значении $l^{(2)}$ величина $l^{(1)}$ определяется из любого из неравенств (\ref{66norm}) не более, чем некоторой константой способов.
 \end{Zam}
Следуя методу Бургейна --- Конторовича, полагается свести оценку мощности множества  $\NN$  к оценке величины $|\M|$. Для этого необходимо выявить  условия, при выполнении которых гарантируется равенство $\NN=\M$.  В прежних работах на эту тему (\cite{FK1} ---  \cite{K5}, \cite{Huang}) проводились построения именно   в таком ключе. В настоящей работе рассматривается 
также ряд случаев, когда равенство $\NN=\M$ не выполняется. 

Прежде всего,   
введем  обозначения
  \begin{equation} 
\label{692f1o6rm}
 \mathbf{P}= \frac{74A^2Q^2_{\alpha}Q_{\beta}}{M_1},{}\ {}\ {}\ 
   \mathbf{T}=\mathbf{T}(\p)=\left [\frac{\mathbf{P}}{\p}\right],{}\ {}\ {}\ 
\M_0=
\sum\limits_{1\le \p\le\mathbf{P}}{}\ 
\sum\limits
 _{0\le{}\ t,T{}\    \le    \mathbf{T}}
|\NN(\p,t,T)|.
\end{equation}

\begin{Le} \label{t8.1} Если выполнено неравенство (\ref{69n1fo2r3m}),
то имеет место оценка
\begin{equation} 
\label{79n4fo2rm1}
|\NN|\le |\M|+\M_0.
\end{equation}
\end{Le} 

 Доказательство. Ввиду первых элементов  минимумов в правых частях неравенств   (\ref{52norm}) и (\ref{52norm1}), имеют место неравенства
 $t,T{}\    \le   \mathbf{T}$. Следовательно, 
 разбивая множество $\NN$ на ряд составляющих его подмножеств, имеем:
\begin{equation} 
\label{69n3fo1rm}
|\NN|\le \sum_{1\le \p \le Q_{\alpha}}{}\ 
\sum_{0\le{}\ t,T{}\    \le   \mathbf{T}}
|\NN(\p,t,T)|.
\end{equation}
Выделяя в (\ref{69n3fo1rm}) слагаемое с   $t=T=0$, приходим к неравенству

$$|\NN|\le 
\sum_{1\le \p \le Q_{\alpha}}{}\ 
|\NN(\p,0,0)|+
\sum_{1\le \p \le Q_{\alpha}}{}\ 
\sum_{0\le{}\ t,T{}\    \le   \mathbf{T}\atop{t^2+T^2\not=0}}
|\NN(\p,t,T)|
.
$$ 
   Но, поскольку выполнены равенства  $\NN(\p,0,0)=\M$ при $\p=q^{(1)}=
q^{(2)}$  и $|\NN(\p,0,0)|=0$ --- в остальных случаях \cite[пп. 6.2.1]{BK} (см. также \cite[доказательство леммы 14.1]{FK3}), то 
\begin{equation} 
\label{69n4fo2rm}
|\NN|\le|\M|
+
\sum_{1\le \p \le Q_{\alpha}}{}\ 
\sum_{0\le{}\ t,T{}\    \le   \mathbf{T}\atop{t^2+T^2\not=0}}
|\NN(\p,t,T)|
.
\end{equation}
Отбросим в (\ref{69n4fo2rm}) равные нулю слагаемые --- те, в которых $\mathbf{T}(\p)<1$.  Так случится при $\p$ из интервала $[[\mathbf{P}]+1, Q_{\alpha}]$. Согласно соотношениям (\ref{69n1fo2r3m}) и  (\ref{692f1o6rm}), выполняется неравенство $\mathbf{P}\le Q_{\alpha}$. Поэтому с помощью приведенных соображений интервал суммирования по $\p$ сокращается, и в результате   получаем оценку (\ref{79n4fo2rm1}). 
Лемма доказана.

\begin{Th}\label{T7.1} Если для любых  натуральных значений $\alpha$ и $\beta$ найдется число $M_1=M_1(\alpha,\beta)$, удовлетворяющее оценкам  (\ref{69n1fo2r3m}) и   
\begin{equation} 
\label{70hokrm11}
\frac{|\M|+\M_0}{|Z||\Omega|^2}
\ll 
\left(M_1
 \right)^{-2+2\delta-\mathbf{c}+O(\varepsilon_0)},
 \end{equation}
то для алфавита $\A$ имеют место формулы (\ref{lc})~---~(\ref{cdd}). 
\end{Th}

Доказательство. С помощью оценок (\ref{79n4fo2rm1}) и  (\ref{70hokrm11}) легко получить неравенство  (\ref{70nokrm}). Тогда утверждение настоящей теоремы следует из теоремы \ref{t7.1}. Теорема доказана. 

\label{8/1}

\section{Оценка величины $|\M|$ }

В этом параграфе будет получена оценка мощности множества $\M$, то есть первого из слагаемых в числителе левой части неравенства (\ref{70hokrm11}).
Напомним, что в этой ситуации выполнено равенство $t=T=0$. В этом случае положим
\begin{equation} 
\label{30norm} 
M_2=\sqrt{\max\left\{1,\frac{1}{4}Q_{\beta-1}\right\}}, {}\ {}\ {}\ {}\ 
M_4=\sqrt{\max\left\{1,\frac{1}{4}Q_{\alpha-1}\right\}}.
 \end{equation}
Далее, для произвольного множества матриц $\Xi\subseteq \Gamma_{\A}$  число решений сравнения 
\begin{equation} 
\label{105norm4}
Uv\equiv uV\pmod q
\end{equation}
в переменных $\left(u\atop{U} \right),\left(v\atop{V} \right)
\in \widetilde{\Xi}$  обозначим через $R_{ q}\left({\Xi}\right)$. Здесь и далее символ $q$ употребляется в том же смысле, что и в
(\ref{17norm}).
   
\begin{Le}\label{9.1} Для всякого числа $M_1$  из интервала 
    (\ref{69n1fo2r3m})
 имеет место неравенство
   \begin{equation} 
\label{109ho11rm}
 |\M|\ll
\left   |\Omega_2\right|
\sum\limits_{\Theta
\in Z}
R_{q}\left(\Omega_3{\Omega}_4\right).
\end{equation}
\end{Le}
Доказательство. Рассмотрим разложения (\ref{1norm}) и (\ref{1niorm}), соответствующие значениям параметров $M_1$, $M_2$ и $M_4$. Неравенством (\ref{69n1fo2r3m}) обеспечивается выполнение оценки (\ref{1yorm1}) для параметров (\ref{30norm}). В \cite[лемма 9.1]{K5} 
 для элементов некоторого аналогичного $\M$ множества было доказано равенство
$g^{(1)}_2=g^{(2)}_2.
$ Тем же методом  осуществляется  доказательство этого равенства  и в случае множества $\M$ для настоящей леммы.

  Далее с минимальными изменениями следует повторить аргументы из доказательства \cite[леммы 9.2] {K5}. Действительно, согласно замечанию \ref{z9.1}, имеют место сравнения  в (\ref{67norm}). Полагая  $g^{(1)}_2=g^{(2)}_2=g_2$, их можно записать в виде 
$$\left(g_2a^{(1)}g_3^{(1)}\widetilde{g}^{(1)}_4\equiv g_2a^{(2)}g_3^{(2)}\widetilde{g}^{(2)}_4 \right)_{1,2}\pmod\q,$$
где индексы ``1,2'' внизу означают  выполнение сравнения по обеим координатам. Отсюда, в виду равенства $\det g_2=1$, получаем: 
\begin{equation} 
\label{103norm}
\left(a^{(1)}g_3^{(1)}\widetilde{g}^{(1)}_4\equiv a^{(2)}g_3^{(2)}\widetilde{g}^{(2)}_4 \right)_{1,2}\pmod\q.
\end{equation}
Положим 
\begin{equation} 
\label{93n1or9m}
g_3^{(1)}\widetilde{g}^{(1)}_4=\left(u \atop{U}\right),{}\ {}\ g_3^{(2)}\widetilde{g}^{(2)}_4=\left(v \atop{V}\right), 
\end{equation}
тогда сравнения (\ref{103norm}) перепишутся в виде
\begin{equation} 
\label{93n1or91m}
a^{(1)}u\equiv a^{(2)}v \pmod\q,{}\ {}\ {}\ {}\ {}\ 
a^{(1)}U\equiv a^{(2)}V \pmod\q.
\end{equation}
Отсюда следует цепочка сравнений:
\begin{equation} 
\label{102norm}
\left( a^{(1)}U\right)v\equiv  a^{(2)}Vv= 
\left( a^{(2)}v\right)V\equiv \left( a^{(1)}u\right)V\pmod\q.
\end{equation}
Но числа $a^{(1)}$ и $\q$ взаимно просты как числитель и знаменатель цепной дроби. Следовательно, сокращая сравнение (\ref{102norm}) на $a^{(1)}$, получаем сравнение (\ref{105norm4}) при $\Xi=\Omega_3\Omega_4$.

  Подытожим сказанное, пересчитывая количество входящих в $\M$ элементов. Для этого выберем и фиксируем матрицу ${g}^{(1)}_2={g}^{(2)}_2$, определяющую элемент $\widetilde{g}^{(1)}_2=\widetilde{g}^{(2)}_2$, одним из $\left   |\Omega_2\right|$ способов --- это первый множитель в (\ref{109ho11rm}). Выберем каким-либо способом   
  число $\Theta=\Theta^{(2)} \in Z,$ представленное в виде  (\ref{50norm}) --- это переменная суммирования в (\ref{109ho11rm}). Тем самым определено число $q=q^{(1)}=q^{(2)}=\q,$
  а также числа $a^{(2)}$ и $l^{(2)}$. Выберем также элементы 
  $$ g^{(1)}_3, {}\ g^{(2)}_3, {}\ g^{(1)}_4, {}\ g^{(2)}_4, {}\ $$
 которые в обозначениях (\ref{93n1or9m}) удовлетворяют сравнению (\ref{105norm4}), одним из 
  $R_{q}\left(\Omega_3{\Omega}_4\right)$ способов.  Заметим теперь, что согласно \cite[доказательству леммы 3.13]{FK3}, число $a^{(1)}$ определяется по $a^{(2)}$  однозначно, исходя из сравнений (\ref{93n1or91m}). Наконец, согласно замечанию \ref{z9.1},  величина $l^{(1)}$ определяется по $l^{(2)}$  из неравенства  (\ref{66norm}) не более, чем некоторой константой способов.
Лемма доказана.

Для каждой  матрицы $g_3\in \Omega_3$ через $g_3{\Omega}_4$ обозначим множество матриц, получающихся умножением $g_3$ на произвольные матрицы из ${\Omega}_4$.

 \begin{Le}\label{9.2}  Для всякого числа $M_1$  из интервала 
    (\ref{69n1fo2r3m})
 имеют место неравенства
  \begin{equation}
\label{108n4orm}
 R_{q}\left(\Omega_3{\Omega}_4\right)\le
 \left|\Omega_3\right|
\sum\limits_{g_3
\in \Omega_3}
R_{q}\left(g_3{\Omega}_4\right)\le
 \left|\Omega_3\right|^2
 \left|\Omega_4\right|
 .
\end{equation}
\end{Le}
Доказательство. Используя обозначения (\ref{93n1or9m}), положим $r=\gcd(v,q)$. Заметим, что числа $v$ и $V$ взаимно просты как числитель и знаменатель цепной дроби. Поэтому, ввиду   сравнения (\ref{105norm4}), выполнено также равенство $r=\gcd(u,q)$. Далее, положим 
$$u_0=\frac{u}{r},{}\ v_0=\frac{v}{r},{}\ 
q_0=\frac{q}{r}.
$$
Введем также обозначения $\left(u_0\right)^{(-1)}$, $\left(v_0\right)^{(-1)}$
 для вычетов по модулю $q_0$, обратных к   $u_0$, $v_0$, соответственно. 
Используя формулу (\ref{11hy}), получаем 
 цепочку равенств 
 
\begin{equation}
\label{108n5orm}
\begin{array}{ll}
R_{q}\left(\Omega_3{\Omega}_4\right)= 
\sum\limits_{\left(u\atop{U} \right),\left(v\atop{V} \right)
\in \widetilde{\Omega}_{3,4}}
  \delta_{q}(Uv- uV) =
\sum\limits_{r \bigl| q}{}\ 
\sum\limits_{\left(u\atop{U} \right),\left(v\atop{V} \right)
\in \widetilde{\Omega}_{3,4}\atop{\gcd(u,q)=\gcd(v,q)=r
}}
\delta_{q_0}
(Uv_0- u_0V)
=\\
\sum\limits_{r \bigl| q}{}\ 
\sum\limits_{\left(u\atop{U} \right),\left(v\atop{V} \right)
\in \widetilde{\Omega}_{3,4}\atop{\gcd(u,q)=\gcd(v,q)=r
}}
\delta_{q_0}
(U\left(u_0\right)^{(-1)}- V\left(v_0\right)^{(-1)})
=\sum\limits_{r \bigl| q}{}\  
\sum\limits^{q_0}_{k=1}
\frac{1}{q_0}\left|
 \sum\limits_{\left(u\atop{U} \right)
\in \Omega_3\widetilde{\Omega}_4\atop{\gcd(u,q)=r
}}
e_{q_0}\left(U\left(u_0\right)^{(-1)}k \right) 
 \right|^2, 
 \end{array}
\end{equation}
где сумма по $r$ берется по всем делителям числа $q$.

Представим аргумент внутренней суммы в (\ref{108n5orm}), стоящей под знаком модуля,  в виде $ \left(u\atop{U} \right)=g_3\widetilde{g}_4,$
где $g_3\in \Omega_3, {}\  {g}_4\in {\Omega}_4,$ и определим функцию $f_{k,r}$ равенством
$$ 
 f_{k,r}(g_3\widetilde{g}_4)
 =e_{q_0}\left(U\left(u_0\right)^{(-1)}k \right) .$$  
 Оценим квадрат модуля внутренней суммы по $ \left(u\atop{U} \right)$ из (\ref{108n5orm}) с помощью неравенств треугольника и Коши ---  Буняковского: 
$$\left|
 \sum\limits_{g_3
\in \Omega_3}
 \sum\limits_{{g}_4
\in {\Omega}_4\atop{\gcd(u,q)=r
}}
f_{k,r}(g_3\widetilde{g}_4) \right|^2
\le
\left(
 \sum\limits_{g_3
\in \Omega_3}
\left| \sum\limits_{{g}_4
\in {\Omega}_4\atop{\gcd(u,q)=r
}}
f_{k,r}(g_3\widetilde{g}_4) \right|\right)^2
\le
 \left|\Omega_3\right|
 \sum\limits_{g_3\in \Omega_3}
\left| \sum\limits_{{g}_4
\in {\Omega}_4\atop{\gcd(u,q)=r
}}
f_{k,r}(g_3\widetilde{g}_4) \right|^2. 
$$
Подставим результат этой оценки в равенство (\ref{108n5orm}):
\begin{equation}
\label{108n5rm7}
R_{q}\left(\Omega_3{\Omega}_4\right) 
\le\left|\Omega_3\right|
 \sum\limits_{r \bigl| q} {}\ 
\sum\limits_{g_3
\in \Omega_3}
\frac{1}{q_0}
 \sum\limits^{q_0}_{k=1}\left|
 \sum\limits_{\left(u\atop{U} \right)
\in g_3\widetilde{\Omega}_4\atop{\gcd(u,q)=r
}}
e_{q_0}\left(U\left(u_0\right)^{(-1)}k \right) 
\right|^2. 
\end{equation}
Раскрывая в   (\ref{108n5rm7}) квадрат модуля тригонометрической суммы и заменяя  сумму по $k$  $\delta$-символом Коробова, получаем:
$$\begin{array}{ll}
R_{q}\left(\Omega_3{\Omega}_4\right)
\le
\left|\Omega_3\right|
 \sum\limits_{r \bigl| q} {}\ 
\sum\limits_{g_3
\in \Omega_3}{}\ 
 \sum\limits_{\left(u\atop{U} \right),\left(v\atop{V} \right)
\in g_3\widetilde{\Omega}_4\atop{\gcd(u,q)=\gcd(v,q)=r
}}
\delta_{q_0}
(Uv_0- u_0V)
=\\{}\\
=\left|\Omega_3\right|
\sum\limits_{g_3
\in \Omega_3}{}\ 
 \sum\limits_{\left(u\atop{U} \right),\left(v\atop{V} \right)
\in g_3\widetilde{\Omega}_4}
\delta_{q}
(Uv- uV)
=
 \left|\Omega_3\right|
\sum\limits_{g_3
\in \Omega_3}
 R_{q}\left(g_3{\Omega}_4\right).
\end{array}
$$
Первое из неравенств в (\ref{108n4orm}) доказано.

Далее, почти дословным повторением рассуждений из \cite[доказательства леммы 9.2]{K5} выводится равномерная оценка по $g_3$ из $\Omega_3$: 
\begin{equation}
\label{108nor8m}
 R_{q}\left(g_3{\Omega}_4\right)\ll \left   |\Omega_4\right|.
\end{equation}
 Остается подставить оценку (\ref{108nor8m}) в доказанное первое из неравенств в   (\ref{108n4orm}), откуда сразу следует второе из них. Лемма доказана.

Рассмотрим неравенство 
\begin{equation}
 \label{70nokrl1}
  \left({Q}_{\alpha}{Q}_{\beta}\right)^{\delta}
   >>
{\left(M_1
 \right)^{2-2\delta+\mathbf{c}+O(\varepsilon_0)}}.
 \end{equation}

\begin{Th}\label{l9.3}  Если число $M_1$ лежит в интервале 
    (\ref{69n1fo2r3m}) и удовлетворяет неравенству 
    (\ref{70nokrl1}), то
 имеют место оценки
\begin{equation}
\label{109n1o11rm}
\frac{|\M|}{|Z|  \left |\Omega\right|^2}
  \ll
 \frac{1}{\left  (Q_{\alpha}Q_{\beta}\right)^{\delta}}
 \left (M_1\right)^{O(\varepsilon_0)}
 \ll  
\left(M_1
 \right)^{-2+2\delta-\mathbf{c}+O(\varepsilon_0)}.
\end{equation}
\end{Th}

Доказательство. 
Ввиду лемм \ref{9.1} и \ref{9.2}, имеет место неравенство
\begin{equation}
\label{108norm}
|\M|  \ll
 \left   |\Omega_2\right|\left   |\Omega_3\right|^2\left   |\Omega_4\right||Z|=\frac{|\Omega|^2|Z|}{\left  |\Omega_2\right|\left   |\Omega_4\right|}
 .
\end{equation}
Подставляя в неравенство (\ref{108norm}) оценки (\ref{2n1orm1}) и (\ref{2n1orm2}), получаем соотношение 
$$ \frac{|\M|}{|Z|  \left |\Omega\right|^2}
 \ll 
 \frac{\left (M_1M_2M_4\right)^{10\varepsilon_0}}{\left (M_2M_4\right)^{2\delta}}.
$$ 
Подставляя сюда значения  $M_2$ и $M_4$ из (\ref{30norm}),  получаем первую из оценок в  (\ref{109n1o11rm}). Вторая из них получается при подстановке неравенства (\ref{70nokrl1}). Теорема доказана.

\section{Определение и свойства соответственных чисел }

 Рассмотрим  неравенства 
\begin{equation} 
\label{69norm}
\frac{\M_0}{|Z||\Omega|^2}
\ll
\frac{1}{\left  (Q_{\alpha}Q_{\beta}\right)^{\delta}}\left (M_1\right)^{O(\varepsilon_0)},
\end{equation}
 
\begin{equation} 
\label{70nokrm11}
\frac{\M_0}{|Z||\Omega|^2}
\ll 
\left(M_1
 \right)^{-2+2\delta-\mathbf{c}+O(\varepsilon_0)}.
  \end{equation}

\begin{Op} Для любой  пары  натуральныx чисел $\alpha$ и $\beta$ всякое число 
$M_1$ со свойствами (\ref{69n1fo2r3m}), (\ref{70nokrl1}) и  (\ref{70nokrm11}) назовем \textbf{соответственным}. 
\end{Op}

\begin{Le}\label{1lL10.1} Если для числа $M_1$ из интервала (\ref{69n1fo2r3m}) выполнены неравенства (\ref{70nokrl1}) и (\ref{69norm}), то имеет место оценка (\ref{70nokrm11}), то есть, число $M_1$ --- соответственное.

Кроме того, если 
   для любых натуральныx чисел $\alpha$ и $\beta$ найдется соответственное  значение $M_1$, 
то для алфавита $\A$ имеют место формулы (\ref{lc})~---~(\ref{cdd}). 
\end{Le}
Доказательство. Из оценок  (\ref{70nokrl1}) и (\ref{69norm})   неравенство (\ref{70nokrm11}) следует непосредственно, так что первая часть леммы доказана. 

Чтобы доказать вторую часть  леммы, сначала применим теорему \ref{l9.3}. Тогда получим, что неравенства  (\ref{109n1o11rm}) выполнены. Остается  подставить оценки  (\ref{109n1o11rm}) и (\ref{70nokrm11}) в теорему \ref{T7.1}.
Лемма доказана.

\begin{Le}\label{t9.1} Пусть для числа $M_1$ выполненено неравенство (\ref{70nokrl1}) и имеет место  хотя бы один из следующих двух наборов соотношений 
\begin{equation} 
\label{71norm}
M_1=120A^2\sqrt{NQ_{\alpha}Q_{\beta}},{}\ {}\ {}\ {}\ 
(Q_{\alpha})^4(Q_{\beta})^4\ge N,
\end{equation}
\begin{equation} 
\label{71nor5m}
M_1=150A^2\left (Q_{\alpha}\right)^{2}Q_{\beta}
,{}\ {}\ {}\ {}\ 
(Q_{\alpha})^{5}(Q_{\beta})^3\le N^2.
\end{equation}
Тогда выполнены формулы (\ref{69n1fo2r3m}) и
${\M_0}=0$ (откуда следует выполнение неравенства (\ref{70nokrm11})), то есть, число $M_1$ --- соответственное.
\end{Le}

Доказательство. Пусть число $M_1$ определено любым из двух перечисленных способов. Тогда оценка (\ref{69n1fo2r3m}) получается применением неравенств из \cite[замечания 7.2]{K5}:
$$  Q_{\alpha}Q_{\beta}\le Q_3\sqrt{N},{}\ {}\ {}\ 
 \left (Q_{\alpha}\right)^2Q_{\beta}\le 3Q_2N,
 $$ справедливых для всех непустых множеств $P_{\alpha,\beta}$.
А 
 равенство $\M_0=0$ получается дословным повторением доказательств \cite[леммы 8.6]{K5} или \cite[леммы 3.4]{FK3}. Лемма доказана.

Всюду далее считаем, что алфавит $\A$ удовлетворяет неравенству 
\begin{equation} 
\label{113sorm}
\frac{11}{14}<\delta\le \frac{4}{5}.
\end{equation}

\begin{Le}\label{t9.11}
Пусть для   пары  натуральныx чисел $\alpha$ и $\beta$  выполнено неравенство 
\begin{equation} 
\label{59n1okrm}
          N^{1-\delta}
          \le   \left  
              (\left (Q_{\alpha}\right)^{2\delta-1}
              \left  (Q_{\beta}\right)^{2\delta-1}
              \right)^{1-\mathbf{c}+O(\varepsilon_0)} .
 \end{equation}
Тогда  найдется соответственное значение $M_1$. 
\end{Le}

Доказательство. Пусть выполнено неравенство (\ref{59n1okrm}).
Поскольку $\delta<\frac{5}{6},$ то выполнено неравенство
\begin{equation} 
\label{89n1okrm}
 \frac{1-\delta}{2\delta-1}>\frac{1}{4}.
\end{equation}
 Отсюда и из  неравенства  (\ref{59n1okrm}) следует неравенство в (\ref{71norm}).
Следовательно, определяя число $M_1$ равенством в (\ref{71norm}), получаем, что неравенство  (\ref{70nokrl1}) совпадает с неравенством  (\ref{59n1okrm}). В таком случае число   $M_1$ является соответственным значением ввиду леммы \ref{t9.1}.  
Лемма доказана.

\begin{Le}\label{T8.1} 
Для всякой  пары  натуральныx чисел $\alpha$ и $\beta$, таких что
$ \alpha\le 5\beta,  
$   найдется соответственное значение $M_1$. 
\end{Le}

Доказательство. Пусть, для начала,  выполнено неравенство
\begin{equation} 
\label{113norm}
N<(Q_{\alpha})^{\frac{5}{2}}(Q_{\beta})^{\frac{3}{2}}.
\end{equation} 
    Тогда  выполнено неравенство в  (\ref{71norm}).  Определим число $M_1$ равенством в  (\ref{71norm}). Тогда, ввиду леммы \ref{t9.1}, достаточно проверить оценку (\ref{70nokrl1}), сводящуюся к неравенству
$$ {\left(NQ_{\alpha}Q_{\beta}\right)^{1-\delta}}
 \ll 
 \left({Q}_{\alpha}{Q}_{\beta}\right)^{\delta+O(\varepsilon_0)-\mathbf{c}}
.
$$ Для этого, ввиду неравенства (\ref{113norm}), достаточно установить оценку
$$ {\left((Q_{\alpha})^{\frac{5}{2}}(Q_{\beta})^{\frac{3}{2}}Q_{\alpha}Q_{\beta}\right)^{1-\delta}}
 \ll 
 \left({Q}_{\alpha}{Q}_{\beta}\right)^{\delta+O(\varepsilon_0)-\mathbf{c}},
$$ 
или
$$ \left(Q_{\alpha}\right)^
{7-9\delta}
\left(Q_{\beta}\right)^
{5-7\delta}
\ll 
 \left({Q}_{\alpha}{Q}_{\beta}\right)^{O(\varepsilon_0)-\mathbf{c}}
 .
$$ Такая оценка следует из  (\ref{113sorm}), поскольку $\delta>\frac{11}{14}>\frac{7}{9},{}\ {}\ $
  $\delta>\frac{5}{7}$.
 
 Если же неравенство (\ref{113norm}) не выполнено, то, следовательно, имеет место противоположное неравенство 
 $(Q_{\alpha})^{\frac{5}{2}}(Q_{\beta})^{\frac{3}{2}}\le N.$
 В этом случае  число $M_1$ определим равенством в (\ref{71nor5m}). Тогда  применение  леммы \ref{t9.1}  при подстановке такого значения $M_1$ в неравенство (\ref{70nokrl1}) приводит  к достаточности проверки неравенства 
$$ {\left(\left (Q_{\alpha}\right)^{2}Q_{\beta}
 \right)^{2-2\delta}}
 \ll 
\left({Q}_{\alpha}\right)^{\delta}
\left({Q}_{\beta}\right)^{\delta}
 \left({Q}_{\alpha}{Q}_{\beta}\right)^{O(\varepsilon_0)-\mathbf{c}},
 $$ или
$$ \left
(Q_{\alpha}\right)^
{4-5\delta}
\left(Q_{\beta}\right)^
{2-3\delta}
\ll 
 \left({Q}_{\alpha}{Q}_{\beta}\right)^{O(\varepsilon_0)-\mathbf{c}}
.
 $$ Переходя  к логарифмам, получаем, что  достаточно  проверить  неравенство  
$$  \frac{\alpha}{\beta} \le  5  <\frac{3\delta-2}{4-5\delta},
  $$ справедливое, ввиду  (\ref{113sorm}), при  $\alpha\le {5}\beta$.
    Лемма доказана. 
  
    Рассмотрим неравенства   
  \begin{equation} 
\label{113tlrm}
\alpha> 5\beta,  
\end{equation}  
\begin{equation} 
\label{59n1krm}
              \left (Q_{\alpha}\right)^\frac{2\delta-1}{1-\delta}
              \left  (Q_{\beta}\right)^\frac{2\delta-1}{1-\delta}
              \le N^{1-\mathbf{c}+O(\varepsilon_0)}              .
 \end{equation}

 \begin{Th}\label{Tcdd} Если  соответственное  значение $M_1$ найдется
   для любых натуральных чисел $\alpha$ и $\beta$, удовлетворяющих  неравенствам (\ref{113tlrm}) и  (\ref{59n1krm}),
то для алфавита $\A$ имеют место формулы (\ref{lc})~---~(\ref{cdd}).  
 \end{Th}
 
 Доказательство. Согласно второй части леммы \ref{1lL10.1}, достаточно найти соответственные числа $M_1$ для каждой пары натуральныx чисел $\alpha$ и $\beta$. 
 Однако при невыполнении неравенства (\ref{113tlrm}) существование соответственного числа доказано в лемме \ref{T8.1}, а при невыполнении неравенства (\ref{59n1krm}) --- в лемме \ref{t9.11}.
 Теорема доказана.
\label{8}

Всюду далее считаем, что неравенства (\ref{113tlrm}) и (\ref{59n1krm}) выполнены. 
    
\begin{Zam} \label{25gf1krm} Для выполнения  неравенства   (\ref{17narb})  и второй из верхних оценок в  (\ref{69n1fo2r3m}) достаточно потребовать, чтобы выполнялось  неравенство 
\begin{equation} 
\label{5gf1krm}
M_1Q_{\alpha}\le 
\left               (\left (Q_{\alpha}\right)
              ^\frac{2\delta-1}{1-\delta}
              \left  (Q_{\beta}\right)^\frac{2\delta-1}{1-\delta}
              \right)
              ^{1-\mathbf{c}+O(\varepsilon_0)} 
.
\end{equation}
Действительно: для неравенства (\ref{17narb})  это утверждение  сразу следует из сравнения неравенств (\ref{17narb}) и (\ref{59n1krm}).

Теми же соображениями обеспечивается  выполнение второй из верхних оценок в  (\ref{69n1fo2r3m}). Именно, чтобы ее получить, достаточно неравенство $M_1\sqrt{Q_{\alpha}Q_{\beta}}\le
M_1Q_{\alpha}$, справедливое ввиду  (\ref{113tlrm}), продолжить  неравенством  (\ref{5gf1krm})  и применить  оценку (\ref{59n1krm}).
\end{Zam}

\section{Оценка величины $\M_0 $}

В этом параграфе будет получена оценка величины $\M_0 $, то есть числителя   левой части неравенства (\ref{70nokrm11}), при невыполнении условий леммы \ref{t9.1}.
 
 Пусть $m$, $n$, $\p$, $t $ --- данные  целые числа, такие что $mn\p$ не равно нулю и числа $m$ и $n$ взаимно просты. Рассмотрим в целых переменных $x$ и $y$ сравнение
\begin{equation}
\label{109nor1m}
 mx-ny\equiv t \pmod{( mn\p )}.
\end{equation}

\begin{Le}\label{l10.1} Найдутся целые числа $x^{(0)}$,  $y^{(0)}$ и $k$, зависящие только от значений параметров $m$, $n$, $\p$, $t $, такие что  
$0\le k<\p$ и 
для любого решения $(x,y)$ сравнения (\ref{109nor1m}) выполняются сравнения
\begin{equation}
\label{110nor1m}
 x\equiv x^{(0)}+kn \pmod{ n\p },{}\ {}\ y\equiv y^{(0)}+km \pmod{m\p }.
{}
\end{equation}
\end{Le}
Доказательство. 
Пусть $(x,y)$, $(X,Y)$ --- какие-нибудь два решения сравнения (\ref{109nor1m}). Тогда подстановка этих значений в исходное сравнение  при последующем вычитании  результатов этих подстановок дает сравнение  
\begin{equation}
\label{111nor1m}
 m(x-X)\equiv n(y-Y) \pmod{ mn\p}.
{}
\end{equation}
Поэтому, ввиду взаимной простоты чисел $m$ и  $n$, выполняются сравнения 
$ x\equiv X\pmod{n}$, $y\equiv Y \pmod{m}.
$
Следовательно, найдутся числа $k_1$ и $k_2$, такие что  
\begin{equation}
\label{113nor1m}
 x= X+k_1 n ,{}\ y= Y+k_2  m.
{}
\end{equation}
Подстановка значений (\ref{113nor1m}) в сравнение (\ref{111nor1m}) приводит к сравнению $k_1\equiv k_2\equiv k\pmod{\p}$, и сравнения (\ref{110nor1m}) доказаны.  Чтобы прийти к неравенству $0\le k<\p$, остается лишь вместо числа $k$ рассмотреть его остаток от деления на  $\p$. Лемма доказана.

Напомним обозначения (\ref{11norm}) и (\ref{50norm}). Пусть фиксированы значения 
\begin{equation}
\label{114nor2m}
  a^{(1)}, a^{(2)}, q^{(1)},   q^{(2)}, t,T,
 {}
\end{equation}
тогда число решений системы из двух сравнений (\ref{69n2f1orm}) в переменных 
\begin{equation}
\label{114gor2m}
  \left( x\atop{X} \right), \left( y\atop{Y} \right) \in    \Omega 
  \end{equation}
обозначим через $X_{t,T}\left(a^{(1)}, a^{(2)}, q^{(1)},   q^{(2)}\right)
 $.
 
Неравенства (\ref{17narb}) для ближайших  далее трех лемм будем считать  
 выполненными, так что, согласно теореме \ref{t6.1}, 
 имеет место разложение полуансамбля  (\ref{1niorm}) с $\Omega_2$ и $\Omega_4$, соответствующими  параметрам (\ref{17narm}). 
 
 \begin{Le}\label{LL9.2} Имеют место оценки
\begin{equation}
\label{14n2r3m11}
 X_{t,T}\left(a^{(1)}, a^{(2)}, q^{(1)},   q^{(2)}\right)
 \ll 
 \left|\Omega_3\right|^2
 \p^2=
 \frac{\left|\Omega\right|^2\p^2}{\left|\Omega_4\right|^2},
\end{equation}
\begin{equation}
\label{114g1or3m}
\left  |\Omega_4\right |>> \left  (Q_{\alpha}\right)^{2\delta
-2\varepsilon_0}.
  \end{equation}  
\end{Le}

Доказательство. При $\alpha=1$ утверждение леммы легко доказывается, поэтому достаточно рассмотреть случай $\alpha>1$. Заметим также, что, согласно  
теореме \ref{t6.1}, выполнено второе из неравенств (\ref{2n1orm}), ввиду которого  выполняется оценка  (\ref{114g1or3m}).

Рассмотрим сравнения в (\ref{69n2f1orm}),   полагая $n=q^{(1)}_0$, $m=q^{(2)}_0$, тогда мы приходим к двум сравнениям вида (\ref{109nor1m}). Так как   $\p n=q^{(1)}$,  $\p m=q^{(2)}$, то,  согласно лемме \ref{l10.1}, найдутся целые числа $k_1$ и $k_2$ в интервале $[0,{} \ \p-1]$ , такие что 
\begin{equation}
\label{114nor3m}
 \left( xa^{(1)}\atop{Xa^{(1)}}\right)\equiv 
 \left( x^{(0)}+k_1q^{(1)}_0\atop{X^{(0)}+k_2q^{(1)}_0}\right) 
 \pmod{q^{(1)}},{}\ {}\ {}\ 
 \left( ya^{(2)}\atop{Ya^{(2)}}\right)\equiv 
 \left( y^{(0)}+k_1q^{(2)}_0\atop{Y^{(0)}+k_2q^{(2)}_0}\right) 
 \pmod{q^{(2)}},
 \end{equation}
 где $x^{(0)}, X^{(0)}$ и $y^{(0)}, Y^{(0)}$ --- константы, аналогичные величинам  $x^{(0)}$  и $y^{(0)}$ из  леммы \ref{l10.1}.

Поскольку $M_2=1$, то $\Omega_2=\{E\}.$ Поэтому вектора $\left(x\atop{X} \right), \left(y\atop{Y} \right)$, участвующие в левых частях сравнений в    (\ref{114nor3m}), можно представить 
в виде
$$ \left(x\atop{X} \right)=g^{(1)}_3\left(x_4\atop{X_4} \right),{}\ 
  \left(y\atop{Y} \right)=g^{(2)}_3\left(y_4\atop{Y_4} \right).
$$
Пусть $\left(a^{(1)}\right)^{-1}$, $\left(a^{(2)}\right)^{-1}$ --- вычеты по модулям $q^{(1)}$ и $q^{(2)}$, обратные к $a^{(1)}$ и $a^{(2)}$, соответственно.
Тогда, умножая сравнения (\ref{114nor3m}) на $\left(a^{(1)}\right)^{-1}$ или, соответственно,  на $\left(a^{(2)}\right)^{-1}$, а также  на матрицы, обратные к  матрицам  $g^{(1)}_3$ или, соответственно, $g^{(2)}_3$, получаем: 
\begin{equation}
\label{114nor4m}
 \left(x_4\atop{X_4} \right)\equiv \left( g^{(1)}_3\right)^{-1}
 \left( x^{(0)}+k_1q^{(1)}_0\atop{X^{(0)}+k_2q^{(1)}_0}\right) 
 \left (a^{(1)}\right)^{-1}\pmod{q^{(1)}},{}\ \\{}\\ 
 \end{equation}
 \begin{equation}
\label{114nor4m1}
\left(y_4\atop{Y_4} \right)\equiv \left( g^{(2)}_3\right)^{-1}
 \left( y^{(0)}+k_1q^{(2)}_0\atop{Y^{(0)}+k_2q^{(2)}_0}\right) 
\left (a^{(2)}\right)^{-1}\pmod{q^{(2)}}.
 \end{equation}
Заметим, что правые части сравнений (\ref{114nor4m}) и  (\ref{114nor4m1}) зависят только от величинт $k_1$ и $k_2$ и значений параметров (\ref{114nor2m}). Поэтому, ввиду сравнений (\ref{114nor4m}) и (\ref{114nor4m1}), при фиксированных значениях $k_1$ и $k_2$ вектора  $\left(x_4\atop{X_4} \right)$ и $\left(y_4\atop{Y_4} \right)$ определены не более, чем однозначно. Это следует из неравенства (\ref{normyy1}), поскольку, согласно выбору числа $M_4$, выполняется неравенство
  ${M_4} < Q_{\alpha-1}\le \min \left\{q^{(1)},q^{(2)}  \right\}.$

 Подытожим сказанное, пересчитывая количество решений полученной системы сравнений. Для этого  матрицы $g^{(1)}_3$ и $g^{(2)}_3$ выберем одним из 
   $\left|\Omega_3\right|^2$ способов, а
         числа $k_1$ и $k_2$ --- одним из
        $\p^2$ вариантов. 
        Поэтому имеет место неравенство (\ref{14n2r3m11}). Лемма доказана.

\begin{Le}\label{lkdL10.1} Пусть 
   для  натуральных  $\alpha$ и $\beta$, для любого значения $\p\in[1,\mathbf{P}]$  и для любых чисел $t,T$ из интервала $[0,\mathbf{T}]$ найдется число $M_1$ со свойствами  (\ref{69n1fo2r3m})  и (\ref{70nokrl1}),  такое что при этом для него выполнено хотя бы одно из неравенств 
\begin{equation} 
\label{70dokrm11}
\frac{|\NN(\p,t,T)|}{\p|Z||\Omega|^2}
\ll 
\frac{\left(M_1
 \right)^{2\delta-\mathbf{c}+O(\varepsilon_0)}}
 {\left({Q}_{\alpha}\right)^{4}  
   \left({Q}_{\beta}\right)^{2}},
  \end{equation}
\begin{equation} 
\label{70cokrm11}
\sum\limits
 _{0\le{}\ t,T{}\    \le    \mathbf{T}\atop{t^2+T^2\not=0}}
\frac{|\NN(\p,t,T)|}{|Z||\Omega|^2}
\ll 
\frac{1}{\p}\left(M_1 \right)^
{2\delta-2 -\mathbf{c}+O(\varepsilon_0)}.
  \end{equation}
Тогда имеет место оценка (\ref{70nokrm11}), то есть, число $M_1$ --- соответственное.
\end{Le}
Доказательство.  Суммируя неравенство (\ref{70dokrm11}) по $t$ и $T$ в пределах от  $0$  до   $\mathbf{T}$, получаем оценку (\ref{70cokrm11}). Далее, суммируя неравенство (\ref{70cokrm11})  по $\p$ в пределах от  $1$  до   $\mathbf{P},$  получаем оценку (\ref{70nokrm11}).
 Лемма доказана.

  \begin{Le} \label{l10.4} Пусть  выполнена оценка (\ref{69n1fo2r3m}). Тогда имеют место оценки
   \begin{equation}
\label{114nor1m11}
\frac{|\NN(\p,t,T)|}{\p|Z|  \left |\Omega\right|^2}
\ll
 \frac{|Z|\p}{\left(Q_{\alpha}\right)^{ 4\delta-4\varepsilon_0}}
\ll
 \frac{ Q_{\beta}}{\left(Q_{\alpha}\right)^{
 4\delta-2 -4\varepsilon_0}}. 
 {}
\end{equation}
\end{Le}

Доказательство. 
    Пересчитаем количество элементов  в $\mathfrak{N} (\p, t,T)$. Для этого    выберем   $|Z|^2$ способами
  числа $\Theta^{(1)},\Theta^{(2)} \in Z$. Тогда значения параметров 
  (\ref{114nor2m}) полностью определены. Следовательно, выбирая $X_{t,T}\left(a^{(1)}, a^{(2)}, q^{(1)},   q^{(2)}\right)
 $   способами
   решение системы из двух сравнений (\ref{69n2f1orm})
в переменных (\ref{114gor2m}), ввиду леммы \ref{LL9.2}, получаем:
$$ |\mathfrak{N} (\p, t,T)|\ll 
 |Z|^2
 X_{t,T}\left(a^{(1)}, a^{(2)}, q^{(1)},   q^{(2)}\right)\le
  \frac{|Z|^2\left|\Omega\right|^2}{\left|\Omega_4\right|^2}
 \p^2.
$$
Применяя здесь неравенство (\ref{114g1or3m}), приходим к первой из оценок в  
(\ref{114nor1m11}). 

Поскольку для чисел $\Theta\in Z$ выполнены соотношения (\ref{17norm}) и (\ref{43norm}), то, следовательно, 
 \begin{equation}
\label{14s1or3m}
  |Z|\le \frac{1}{\p}\left( Q_{\alpha} \right)^2Q_{\beta}.
   \end{equation} 
   Это неравенство получается  из учета не более $ \frac{1}{\p}\left( Q_{\alpha} \right)^2$ дробей $\frac{a}{\q}$, в которых $\q$ делится на $\p$, и не более $ Q_{\beta}$ значений параметра $l$ --- для каждой из них. Подставляя оценку (\ref{14s1or3m}) в доказанную первую из оценок в (\ref{114nor1m11}), получаем вторую из них. Лемма доказана. 

Напомним, что здесь и далее  неравенства (\ref{113sorm}), (\ref{113tlrm}) и (\ref{59n1krm}) выполнены.

\begin{Th} \label{1g4nort} Пусть выполнено неравенство 
\begin{equation}
\label{1g4or1m}
 \delta> \frac{\sqrt{21}-3}{2}=0.791\ldots   .
\end{equation}
Тогда   для алфавита $\A$ имеют место формулы (\ref{lc})~---~(\ref{cdd}). 
 \end{Th}

Доказательство. В силу теоремы \ref{Tcdd} и леммы \ref{lkdL10.1}, достаточно доказать соответственность  числа $M_1$, заданного равенством 
\begin{equation}
\label{1g4nor1m}
M_1=
\left(
\left(Q_{\alpha}\right)^{\frac{3-2\delta}{\delta}}
\left(Q_{\beta}\right)^{\frac{3}{2\delta}}
\right)^{1+\mathbf{c}+O(\varepsilon_0)}.
\end{equation}
То есть, достаточно доказать неравенства (\ref{17narb}), (\ref{69n1fo2r3m}),   (\ref{70nokrl1}) и (\ref{70dokrm11}) для такого $M_1$.

   Для этого рассмотрим неравенство
\begin{equation}
\label{11vvo2r72m}
\left(Q_{\alpha}\right)^
{\frac{3-\delta}{\delta}}
\left(Q_{\beta}\right)^{\frac{3}{2\delta}}
\le  
\left  
              (\left (Q_{\alpha}\right)
              ^\frac{2\delta-1}{1-\delta}
              \left  (Q_{\beta}\right)^\frac{2\delta-1}{1-\delta}
              \right)^{1-\mathbf{c}+O(\varepsilon_0)} 
          . 
 \end{equation}
 Оно доказывается с помощью  неравенств
\begin{equation}
\label{11veo2r72m}
\frac{3-\delta}{\delta}<\frac{2\delta-1}{1-\delta},{}\ {}\ {}\ {}\ {}\ {}\ 
\frac{3}{2\delta}<\frac{2\delta-1}{1-\delta},
\end{equation}
имеющих место ввиду (\ref{1g4or1m}) (действительно, неравенства (\ref{11veo2r72m}) можно преобразовать, соответственно, к виду $\delta^2+3\delta-3>0$ или $4\delta^2+\delta-3>0$).  Отсюда получается неравенство (\ref{5gf1krm}). Поэтому, согласно  замечанию  \ref{25gf1krm}, имеют место неравенство (\ref{17narb}) и вторая из верхних оценок в  (\ref{69n1fo2r3m}). Остальные оценки в (\ref{69n1fo2r3m}) доказываются применением неравенств
$$1<\frac{3-2\delta}{\delta}<5,{}\ {}\ {}\ {}\ {}\ {}\ 
1<\frac{3}{2\delta}<5,
$$
имеющих место ввиду оценки (\ref{113sorm}).
   
      Следовательно, выполнены условия леммы \ref{l10.4}, согласно которой  выполняется  оценка (\ref{114nor1m11}). Чтобы обеспечить выполнение условий леммы \ref{lkdL10.1}, нужно, в частности, вывести из доказанной оценки (\ref{114nor1m11}) неравенство (\ref{70dokrm11}).
 Однако неравенства  (\ref{70dokrm11})  и (\ref{114nor1m11}) соответствуют одно другому при $M_1,$ определенном  равенством  (\ref{1g4nor1m}).
 
 Рассмотрим следующее из условия  (\ref{1g4or1m}) неравенство
\begin{equation}
\label{g4nor1m}
\left(Q_{\alpha}\right)^{\frac{3-2\delta}{\delta}}
\left(Q_{\beta}\right)^{\frac{3}{2\delta}}
\le
\left(
\left(Q_{\alpha}\right)^{\frac{\delta}{2-2\delta}}
\left(Q_{\beta}\right)^{\frac{\delta}{2-2\delta}}
\right)^{1-\mathbf{c}+O(\varepsilon_0)}.
\end{equation}
(Действительно, для его доказательства достаточно лишь проверить неравенства 
 $$  \frac{3-2\delta}{\delta}<\frac{\delta}{2-2\delta},{}\ {}\ {}\ {}\ {}\ {}\ 
\frac{3}{2\delta}<\frac{\delta}{2-2\delta},
                                                 $$  
                                                 равносильные, соответственно,  оценкам  
    $\delta^2-10\delta+6<0$ или $3\delta^2+3\delta-3>0$, справедливым при $\delta>\frac{5-\sqrt{7}}{3}=0.7847\ldots$ или, соответственно, при выполнении неравенства (\ref{1g4or1m}) в точности.)                                         
                                                 Из
   неравенства (\ref{g4nor1m}) оценка (\ref{70nokrl1})  
   следует непосредственно. Все условия леммы \ref{lkdL10.1} выполнены.
   Поэтому утверждение теоремы следует из леммы \ref{lkdL10.1}.
Теорема доказана.

\section{Обобщение теоремы \ref{1g4nort}.}
Обобщим результат этой теоремы на величины $\delta$, не удовлетворяющие неравенству (\ref{1g4or1m}).
 Для этого рассмотрим действительные  параметры 
$\mathbf{m}$ и  $\rho=\rho(\mathbf{m})$ и  соотношения \begin{equation}
 \label{90knokrl19}
  0<2\rho<\delta,
 \end{equation}
\begin{equation}
 \label{70nokrl19}
  \left({Q}_{\alpha}{Q}_{\beta}\right)^{\rho}
   >>
{\left(M_1
 \right)^{1-\delta+\mathbf{c}+O(\varepsilon_0)}},
 \end{equation}  
\begin{equation} 
\label{69norm19}
\frac{\M_0}{|Z||\Omega|^2}
\ll
\frac{1}{\left  (Q_{\alpha}Q_{\beta}\right)^{2\rho}}\left (M_1\right)^{O(\varepsilon_0)},
\end{equation}

 \begin{equation} 
\label{69bgrm19}
\frac{\alpha}{\beta}= \mathbf{m}> 5.
\end{equation}

\begin{Op} Пусть действительное число $\mathbf{m}> 5$ и натуральные числа $\alpha$ и $\beta$ подчинены условию (\ref{69bgrm19}), а для  числа $M_1$ из интервала (\ref{69n1fo2r3m}) выполнена оценка (\ref{70nokrm11}). Если найдется число $2\rho$ из интервала (\ref{90knokrl19}), такое что выполнено неравенство (\ref{70nokrl19}),  то число 
$M_1$  назовем 
$\mathbf{m}$- \textbf{соответственным}. 
\end{Op}

\begin{Le}\label{l9.31} Пусть для натуральных чисел $(\alpha,\beta)$, подчиненных условию  (\ref{69bgrm19}), найдется $\mathbf{m}$-соответственное значение $M_1$. Тогда
  имеют место оценки (\ref{109n1o11rm}).
\end{Le}

Доказательство.   Возводя   неравенство  
 (\ref{70nokrl19}) в квадрат и умножая на неравенство $\left(Q_{\alpha}Q_{\beta}\right)^{\delta-2\rho}
>  1 
 ,$
   выполненное по условию (\ref{90knokrl19}), получаем 
 оценку (\ref{70nokrl1}). 
 Так что утверждение  леммы следует из теоремы  \ref{l9.3} непосредственно. Лемма доказана.

\begin{Le}\label{lL10.1} Пусть фиксировано  число $\mathbf{m}> 5$. Если для числа $M_1$ из интервала (\ref{69n1fo2r3m}) выполнены неравенства (\ref{70nokrl19}) и (\ref{69norm19}) с некоторым подчиненным условию (\ref{90knokrl19})   числом $\rho$, то имеет место оценка (\ref{70nokrm11}), то есть число $M_1$ --- $\mathbf{m}$-соответственное (для любой  пары натуральных чисел $\alpha$ и $\beta$, подчиненных условию (\ref{69bgrm19})).

Кроме того, если для любого числа $\mathbf{m}> 5$ и
   для любых натуральных чисел $\alpha$ и $\beta$, подчиненных условию  (\ref{69bgrm19}), найдется     $\mathbf{m}$-соответственное  значение $M_1$, 
то для алфавита $\A$ имеют место формулы (\ref{lc})~---~(\ref{cdd}). 
\end{Le}

Доказательство. Возводя оценку (\ref{70nokrl19}) в $(-2)-$ю степень и подставляя в правую часть (\ref{69norm19}), получаем   неравенство (\ref{70nokrm11}). Так что первая часть леммы доказана. 

Переходим ко второй части. Согласно лемме \ref{l9.31},
 неравенства  (\ref{109n1o11rm}) выполнены. Остается подставить оценки  (\ref{109n1o11rm}) и (\ref{70nokrm11}) в теорему \ref{T7.1}.
Лемма доказана.

\begin{Le}\label{t90.1} Если для числа $M_1$ имеют место   соотношения  (\ref{71norm}) и (\ref{70nokrl19}),
то выполнено  равенство ${\M_0}=0$ (откуда следует неравенство (\ref{69norm19}) с любым подчиненным условию (\ref{90knokrl19})   числом $\rho$) и  неравенство (\ref{69n1fo2r3m}), то есть, ввиду первой части леммы \ref{lL10.1}, число $M_1$ --- $\mathbf{m}$-соответственное.
\end{Le}

Доказательство настоящей леммы состоит в почти дословном повторении доказательства леммы \ref{t9.1}. Лемма доказана.

\begin{Le}\label{t9.1h1}
Пусть для     натуральныx чисел $\alpha$ и $\beta$, подчиненных условию (\ref{69bgrm19}),  выполнены   неравенства
\begin{equation} 
\label{59n1okrm65}
          N^{1-\delta}
          \le   \left  
              (\left (Q_{\alpha}\right)^{2\rho+\delta-1}
              \left  (Q_{\beta}\right)^{2\rho+\delta-1}
              \right)^{1-\mathbf{c}+O(\varepsilon_0)} ,
 \end{equation}
\begin{equation} 
\label{59h1okrm}
\frac{2\rho+\delta-1}{1-\delta}<4  .
               \end{equation}
Тогда  найдется $\mathbf{m}$-соответственное значение $M_1$. 
\end{Le}

Доказательство. Из  неравенств  (\ref{59n1okrm65}) и  (\ref{59h1okrm}) следует неравенство в (\ref{71norm}).
Следовательно, определяя число $M_1$ равенством в (\ref{71norm}), получаем, что из неравенства  (\ref{59n1okrm65}) следует неравенство  (\ref{70nokrl19}). В таком случае число   $M_1$ является $\mathbf{m}$-соответственным значением ввиду леммы \ref{t90.1}.  
Лемма доказана. 

   Рассмотрим оценку  
\begin{equation} 
\label{59nyy1kr}
 \left (Q_{\alpha}\right)^{\frac{2\rho+\delta-1}{1-\delta}}
              \left  (Q_{\beta}\right)^{\frac{2\rho+\delta-1}{1-\delta}}              \le N^{1+O(\varepsilon_0)}              .
 \end{equation}

 \begin{Le}\label{Tcdd} Пусть выполнено неравенство (\ref{59h1okrm}). Если  
$\mathbf{m}$-соответственное  значение $M_1$ найдется
   для любых натуральных чисел $\alpha$ и $\beta$, удовлетворяющих   неравенствам (\ref{69bgrm19}) и  (\ref{59nyy1kr}),
то для алфавита $\A$ имеют место формулы (\ref{lc})~---~(\ref{cdd}).  
 \end{Le}
 
 Доказательство. Согласно второй части леммы \ref{lL10.1}, достаточно найти $\mathbf{m}$-   
 соответственные числа $M_1$ для каждой пары натуральныx чисел $\alpha$ и $\beta$. 
 Однако при невыполнении неравенства (\ref{59nyy1kr}) существование соответственного числа доказано в лемме \ref{t9.1h1}.
 Лемма доказана.
\label{8}

Всюду далее считаем, что неравенства  (\ref{59h1okrm}) и (\ref{59nyy1kr}) выполнены.

\begin{Zam} \label{200f1krm} Для выполнения  неравенства   (\ref{17narb})  и второй из верхних оценок в  (\ref{69n1fo2r3m}) достаточно потребовать, чтобы выполнялось  неравенство 
\begin{equation} 
\label{5g091krm}
M_1Q_{\alpha}\le 
\left (
\left (Q_{\alpha}\right)^{\frac{2\rho+\delta-1}{1-\delta}}
              \left  (Q_{\beta}\right)^{\frac{2\rho+\delta-1}{1-\delta}} \right)^{1-\mathbf{c}+O(\varepsilon_0)}
              .
\end{equation}
Действительно: для неравенства (\ref{17narb})  это утверждение  сразу следует из сравнения неравенств (\ref{17narb}) и (\ref{59nyy1kr}).

Теми же соображениями обеспечивается  выполнение второй из верхних оценок в  (\ref{69n1fo2r3m}). Именно, чтобы ее получить, достаточно неравенство $M_1\sqrt{Q_{\alpha}Q_{\beta}}\le
M_1Q_{\alpha}$, справедливое ввиду  (\ref{69bgrm19}), продолжить неравенством  (\ref{5g091krm}) и применить оценку (\ref{59nyy1kr}).
\end{Zam}

\begin{Le}\label{lkL170.1} Пусть для некоторых натуральных $\alpha$ и $\beta$, таких что выполнено равенство (\ref{69bgrm19}), 
   найдется подчиненное условию (\ref{90knokrl19})  число  $\rho$, для которого существует 
     число $M_1$ со свойствами  (\ref{69n1fo2r3m})  и (\ref{70nokrl19}). Пусть при этом  
   для любого значения $\p\in[1,\mathbf{P}]$  и для любых $t$ и $ T$ из интервала $[0,\mathbf{T}]$  выполнено неравенство 
\begin{equation} 
\label{70dorm11}
\frac{|\NN(\p,t,T)|}{\p|Z||\Omega|^2}
\ll 
\frac{
\left(M_1
 \right)^{2+O(\varepsilon_0)}}
 {\left({Q}_{\alpha}\right)^{4+2\rho}  
   \left({Q}_{\beta}\right)^{2+2\rho}}.
  \end{equation}
 Тогда имеет место оценка (\ref{70nokrm11}), то есть, число $M_1$ --- $\mathbf{m}$-соответственное.
\end{Le}
Доказательство.  Суммируя неравенство (\ref{70dorm11}) по $t$ и $T$ в пределах от  $0$  до   $\mathbf{T}$ и   по $\p$ в пределах от  $1$  до   $\mathbf{P},$  получаем оценку (\ref{69norm19}). Поэтому утверждение леммы следует из первой части леммы \ref{lL10.1}. 
 Лемма доказана.

Всюду далее число $M_1$ задается равенством 
\begin{equation}
\label{1g5or1m} 
M_1=\left(
 \left(Q_{\alpha}\right)^{3-{2}\delta+\rho}
\left(Q_{\beta}\right)^{1.5+
\rho}
\right)
^{1+\mathbf{c}+O(\varepsilon_0)}.
\end{equation}

\begin{Le} \label{1g4nort1}   Пусть для некоторого числа $\mathbf{m}> 5$  найдется число  $\rho$, такое что  выполнены неравенства  (\ref{90knokrl19}),  (\ref{59h1okrm}) и еще --- два:
\begin{equation} 
\label{1g8o1m2}
{\mathbf{m}}
 \left({4}\delta-10-{2\rho}+\frac{6}{\delta}\right)
 <   2\rho+3-\frac{3}{\delta}
,
\end{equation}
\begin{equation} 
\label{1g8o762}
 {4}-2\delta+\rho-\frac{2\rho+\delta-1}{1-\delta}
 <0<\frac{2\rho+\delta-1}{1-\delta}
   -1.5-\rho.
\end{equation}

Тогда число $M_1$, заданное равенством  (\ref{1g5or1m}),
   является $\mathbf{m}$-соответственным  
   для любых натуральных $\alpha$ и $\beta$, подчиненных условию  (\ref{69bgrm19}).
\end{Le}
Доказательство. Предполагается использовать леммы \ref{l10.4} и  \ref{lkL170.1}, поэтому
следует обеспечить условия их применимости,  доказав 
 неравенства  (\ref{17narb}), (\ref{69n1fo2r3m}),   (\ref{70nokrl19}) и (\ref{70dorm11}). 
    
   Для начала    заметим, что, ввиду (\ref{1g5or1m}) и (\ref{1g8o762}), выполнена  оценка (\ref{5g091krm}).
 Из нее, согласно  замечанию  \ref {200f1krm}, следуют неравенство (\ref{17narb}) и вторая из верхних оценок в  (\ref{69n1fo2r3m}). 
 Остальные оценки в (\ref{69n1fo2r3m}) доказываются применением следующих из условия (\ref{90knokrl19})  неравенств
$$1<3-{2}\delta+
 {\rho}<5,{}\ {}\ {}\ {}\ {}\ {}\ 
1<1.5 +{\rho}<5.
$$
  Следовательно, выполнены условия леммы \ref{l10.4}. Согласно этой лемме, ввиду оценки (\ref{114nor1m11}), равенством  
(\ref{1g5or1m}) обеспечивается выполнение неравенства (\ref{70dorm11}). 

  Далее, заметим, что неравенство 
    (\ref{1g8o1m2}) можно с помощью равенства в (\ref{69bgrm19}) преобразовать к виду   
\begin{equation} 
\label{1g4o1m2}
\alpha \rho+\beta \rho
>\left(\alpha
\left(3-{2}\delta+{\rho}\right)
 +\beta
\left(1.5 +{\rho}\right)\right)
(1-\delta)
.
\end{equation}
Беря экспоненту от обеих частей неравенства (\ref{1g4o1m2}), получаем оценку (\ref{70nokrl19}):
$$\left(Q_{\alpha}Q_{\beta}\right)^{\rho}
>>\left(
\left(Q_{\alpha}\right)^{
3-2\delta+{\rho}}
\left(Q_{\beta}\right)^{
1.5 +{\rho}}
\right)^{1-\delta}
\left(
Q_{\alpha}
Q_{\beta}
\right)^{\mathbf{c}+O(\varepsilon_0)}.
$$  
 Поэтому выполнены условия  леммы \ref{lkL170.1}, согласно которой
число $M_1$ --- $\mathbf{m}$-соответственное. Лемма доказана.

\begin{Th} \label{TH1g4ort1}  Пусть для некоторого числа $\mathbf{m}> 5$  найдется число  $\rho$, такое что  выполнены неравенства  
\begin{equation} 
\label{1g6oj1m1}
2\rho(1+\mathbf{m})>
\max\left\{
\frac{4\delta^2-10\delta+6}{1+\delta}(\mathbf{m}+1)
,{}\ {}\ {}\ 
\frac{5-5\delta}{\delta+1}(\mathbf{m}+1)
,{}\ {}\ {}\ \mathbf{m}\left({4}\delta- 10+\frac{6}{\delta}\right)-3+\frac{3}{\delta}
 \right\}
            ,
\end{equation}
\begin{equation} 
\label{1g8oj1m2}
 2\rho(1+\mathbf{m})
<\min\left\{\delta(\mathbf{m}+1),{}\ {}\ {}\ 
(5-5\delta)(\mathbf{m}+1)
\right\}
 .
\end{equation}
Тогда число $M_1$, заданное равенством  (\ref{1g5or1m}),
   является $\mathbf{m}$-соответственным  
   для любых натуральных $\alpha$ и $\beta$, подчиненных условию  (\ref{69bgrm19}).
\end{Th}
Доказательство. Решая относительно переменной $\rho$ неравенства  (\ref{90knokrl19}),  (\ref{59h1okrm}),
(\ref{1g8o1m2}) и 
(\ref{1g8o762}),
 получаем неравенства
 (\ref{1g6oj1m1}) и (\ref{1g8oj1m2}). Поэтому утверждение леммы следует из леммы \ref{1g4nort1}. 
 Теорема доказана.

\section{Доказательство теорем \ref{2.1} и \ref{2.2}}
Расшифруем неравенство (\ref{1g8oj1m2}). Минимум в правой части этого неравенства достигается на первом элементе ввиду неравенства 
(\ref{113sorm}): действительно, сравнение элементов этого минимума приводит к неравенству $\delta  <     \frac{5}{6}          $. 

Далее, система неравенств (\ref{1g6oj1m1}) и (\ref{1g8oj1m2})  относительно переменной  $2\rho(1+\mathbf{m}) $ разрешима тогда и только тогда, когда нижняя оценка этой переменной, данная в  неравенстве (\ref{1g6oj1m1}), меньше ее верхней оценки, данной  неравенством  (\ref{1g8oj1m2}). Другими словами, когда выполнено неравенство 
\begin{equation} 
\label{1g87j1m1}
\delta(1+\mathbf{m})>
\max\left\{
\frac{4\delta^2-10\delta+6}{1+\delta}(\mathbf{m}+1)
,{}\ {}\ {}\ 
\frac{5-5\delta}{\delta+1}(\mathbf{m}+1)
,{}\ {}\ {}\ \mathbf{m}\left({4}\delta- 10+\frac{6}{\delta}\right)-3+\frac{3}{\delta}
 \right\}
            .
\end{equation}

Покажем, что первые два элемента максимума в (\ref{1g87j1m1}) могут быть отброшены. Для этого достаточно доказать выполнение неравенств
\begin{equation} 
\label{187j1m1}
\delta>
\frac{4\delta^2-10\delta+6}{1+\delta}           ,{}\ {}\ {}\ 
\delta
>
\frac{5-5\delta}{\delta+1}
            .
\end{equation}
Оба неравенства справедливы 
 ввиду оценки 
(\ref{113sorm}): действительно, эти неравенства приводят к соотношениям $3\delta^2-11\delta+6<0$ или, соответственно, $\delta^2+6\delta-5>0$, справедливым при $\delta>\frac{2}{3}$ или, соответственно, при $\delta>\sqrt{14}-{3}=0.741\ldots$. 

Таким образом, для проверки неравенства (\ref{1g87j1m1}) остается установить, что

\begin{equation} 
\label{1g807j1m1}
\delta(1+\mathbf{m})>
\mathbf{m}\left({4}\delta- 10+\frac{6}{\delta}\right)-3+\frac{3}{\delta}
            .
\end{equation}
После упрощения неравенство (\ref{1g807j1m1}) сводится к неравенству 
\begin{equation} 
\label{1807j1m1}
 \mathbf{m}>5> \frac{\delta^2+3\delta-3}
{3\delta^2-10\delta+6},    
  \end{equation}
справедливому ввиду неравенства 
(\ref{113sorm}): действительно,  оценка (\ref{1807j1m1}) приводит к неравенству 
$14\delta^2-53\delta+33<0$, справедливому при $\delta>\frac{11}{14}$.

Следовательно, система неравенств (\ref{1g6oj1m1}) и (\ref{1g8oj1m2})  относительно переменной  $2\rho(1+\mathbf{m}) $ разрешима. 
Остается лишь применить вторую часть леммы \ref{lL10.1} и теорему \ref{TH1g4ort1}. 
Теоремы доказаны.

\end{document}